\documentclass[12pt]{amsart}

\usepackage{amsthm}
\usepackage{amssymb}
\usepackage{amsmath}

\newcommand{\DAHA}{\mathbb{H}}
\newcommand{\dDAHA}{\overline{\mathbb{H}}}
\numberwithin{equation}{section}
\newtheorem{thm}{Theorem}[section]
\newtheorem{prop}[thm]{Proposition}
\newtheorem{lem}[thm]{Lemma}

\theoremstyle{definition}
\newtheorem{definition}[thm]{Definition}
\newtheorem{defn}[thm]{Definition}
\newtheorem{rem}[thm]{Remark}

\begin{document}

\title[Differential equations compatible with boundary rational qKZ]
{Differential equations compatible with \\ boundary rational qKZ equation}

\author{Yoshihiro Takeyama}
\address{Department of Mathematics, 
Graduate School of Pure and Applied Sciences, 
Tsukuba University, Tsukuba, Ibaraki 305-8571, Japan}
\email{takeyama@math.tsukuba.ac.jp}

\dedicatory{\it Dedicated to Professor Tetsuji Miwa on his sixtieth birthday} 

\begin{abstract}
We give differential equations compatible with the rational qKZ equation 
with boundary reflection. 
The total system contains the trigonometric degeneration of the bispectral qKZ equation 
of type $(C_{n}^{\vee}, C_{n})$ which in the case of type $GL_{n}$ was studied by van Meer and Stokman. 
We construct an integral formula for solutions to our compatible system 
in a special case. 
\end{abstract}
\maketitle

\setcounter{section}{0}
\setcounter{equation}{0}


\section{Introduction}

In this paper we give differential equations compatible with 
the rational version of the quantum Knizhnik-Zamolodchikov (qKZ) equation \cite{FR} 
with boundary reflection, which we call the {\it boundary rational qKZ equation}. 

Let $V=\mathbb{C}^{2N}$ be a vector space of even dimension. 
The boundary rational qKZ equation is the following system of difference equations 
for an unknown function $f(x\,|\,y)$ on $(\mathbb{C}^{\times})^{N}\times \mathbb{C}^{n}$ 
taking values in $V^{\otimes n}$:  
\begin{align*}
& 
f(x\,|\ldots , y_{m}-c, \ldots) \\
&=R_{m, m-1}(y_{m}-y_{m-1}-c) \cdots R_{m, 1}(y_{m}-y_{1}-c) 
K_{m}(y_{m}-c/2 \, | \, x, \beta) \\ 
&\times 
R_{1,m}(y_{1}+y_{m}) \cdots R_{m-1,m}(y_{m-1}+y_{m}) 
R_{m,m+1}(y_{m}+y_{m+1}) \cdots R_{m, n}(y_{m}+y_{n}) \\ 
&\times 
K_{m}(y_{m} \, | \, \underbar{$1$}, \alpha)\, 
R_{m, n}(y_{m}-y_{n}) \cdots R_{m, m+1}(y_{m}-y_{m+1}) 
f(x\,|\ldots , y_{m}, \ldots) 
\end{align*} 
for $1 \le m \le n$,  
where $c, \alpha$ and $\beta$ are parameters and 
$\underbar{$1$}=(1, \ldots , 1) \in (\mathbb{C}^{\times})^{N}$. 
The linear operator $R(\lambda)$ on $V^{\otimes 2}$ is 
the rational $R$-matrix, 
and $K(\lambda\,|\,x,\beta) \in {\rm End}(V)$ is the boundary $K$-matrix 
which is a linear sum of the identity and the reflection of the basis of $V$ 
with exponent $x=(x_{1}, \ldots , x_{N})$ 
(see \eqref{eq:def-K} below).  
The indices of $R$ and $K$ in the right hand side signify the position of the components 
of $V^{\otimes n}$ on which they act. 
The boundary rational qKZ equation can be regarded as a slight generalization 
of an additive degeneration of 
Cherednik's trigonometric qKZ equation \cite{C1} 
associated with the root system of type $C_{n}$. 
Such equation was also derived as that for correlation functions of 
the integrable spin chains with a boundary \cite{JKKMW, C1.5}. 
Recently combinatorial aspects of a special polynomial solution 
have attracted attention \cite{Ca, DZ, GPS}. 

In this paper we give commuting differential operators in the form 
\begin{align*}
D_{a}(x\,|\,y)=cx_{a}\frac{\partial}{\partial x_{a}}+L_{a}(x\,|\,y) \qquad (1 \le a \le N), 
\end{align*}
where $L_{a}(x\,|\,y) \, (1 \le a \le N)$ are commuting linear operators acting on $V^{\otimes n}$, 
and prove that the system consisting of the boundary rational qKZ equation and 
the differential equations $D_{a}(x\,|\,y)f(x\,|\,y)=0 \, (1 \le a \le N)$ is compatible
(see Theorem \ref{thm:main} below). 
Such compatible system is obtained for the differential KZ or the qKZ equations without boundary reflection 
by Etingof, Felder, Markov, Tarasov and Varchenko in more general settings \cite{EV, FMTV, TV1, TV2}.  

In \cite{MS} van Meer and Stokman constructed a consistent system of 
$q$-difference equations, which they call the bispectral qKZ equation, 
using the double affine Hecke algebra (DAHA) \cite{C2} of type $GL_{n}$. 
The key ingredients are Cherednik's intertwiners and the dual anti-involution. 
As mentioned in \cite{MS}, their construction can be extended to arbitrary root system. 
In this paper we consider the case of type $(C_{n}^{\vee}, C_{n})$. 
The DAHA of this type also has intertwiners and dual anti-involution \cite{S}, 
and hence we can construct the bispectral qKZ equation. 
Now recall that the DAHA has a trigonometric degeneration \cite{C2}. 
In this degeneration the bispectral qKZ equation turns into a system of 
differential equations called the affine KZ equation (see Section 1.1.3 of \cite{C2}) 
and additive difference equations. 
We prove that the system is contained in our compatible system 
of differential equations and the boundary rational qKZ equation with $N=n$ 
restricted to a subspace of $V^{\otimes n}$ isomorphic to 
the group algebra of the Weyl group of type $C_{n}$. 
As will be seen in Section \ref{subsec:embedding} 
the differential operator $D_{a}(x\,|\,y)$ does not literally appear 
in the affine KZ equation because some of its parts act by zero on the subspace. 
Thus our compatible system gives a non-trivial generalization of the trigonometric degeneration 
of the bispectral qKZ equation of type $(C_{n}^{\vee}, C_{n})$. 


As previously mentioned 
the boundary rational qKZ equation can be regarded as 
an additive degeneration of the qKZ equation on the root system of type $C_{n}$. 
For the qKZ equation of type $C_{n}$, 
Mimachi obtained an integral formula for solutions in a special case \cite{M}. 
Similar construction works for the boundary rational qKZ equation with 
$\alpha=\beta=k/2$, where $k$ is a parameter contained in the $R$-matrix, 
and the exponent $x$ restricted to the hyperplane $x_{2}=\cdots =x_{N}=1$. 
We prove that the solutions obtained in this way satisfy the differential equation 
$D_{1}(x_{1}, 1, \ldots , 1 \, |\, y)f=0$. 
Thus we get solutions of our compatible system in a special case. 

The rest of this paper is organized as follows. 
In Section \ref{sec:bdry-qkz} we give the definition of the boundary rational qKZ equation. 
In Section \ref{sec:diff-eqn} we construct the commuting differential operators and prove that 
they are compatible with the boundary rational qKZ equation. 
In Section \ref{sec:bispectral} we calculate the trigonometric degeneration 
of the bispectral qKZ equation of type $(C_{n}^{\vee}, C_{n})$ 
and see that it is contained in our compatible system. 
In Section \ref{sec:integral-formula} we give an integral formula for solutions 
to our compatible system with $\alpha=\beta=k/2$ and $x=(x_{1}, 1, \ldots , 1)$.


\section{Boundary rational qKZ equation}\label{sec:bdry-qkz}

Let $n$ and $N$ be positive integers, and  
$c, k, \alpha$ and $\beta$ be generic non-zero complex numbers. 
Let $V=\oplus_{a=1}^{N}\left(\mathbb{C}v_{a} \oplus \mathbb{C}v_{\overline{a}}\right)$ 
be a vector space with the basis 
$\{v_{l}\}_{ l=1, \ldots , N, \bar{1}, \ldots , \bar{N}}$. 

The rational $R$-matrix acting on $V^{\otimes 2}$ is defined by 
\begin{align*}
R(\lambda):=\frac{\lambda+kP}{\lambda+k}, 
\end{align*}
where $P$ is the transposition $P(u \otimes v):=v \otimes u$. 
It is a rational solution of the Yang-Baxter equation 
\begin{align}
& \label{eq:YBE} 
R_{12}(\lambda_{1}-\lambda_{2})R_{13}(\lambda_{1}-\lambda_{3})R_{23}(\lambda_{2}-\lambda_{3}) \\ 
& \quad {}=
R_{23}(\lambda_{2}-\lambda_{3})R_{13}(\lambda_{1}-\lambda_{3})R_{12}(\lambda_{1}-\lambda_{2}) 
\nonumber
\end{align}
on $V^{\otimes 3}$, 
where $R_{ij}(\lambda)$ is the linear operator acting on the tensor product of the $i$-th and the $j$-th 
component of $V^{\otimes 3}$ as $R(\lambda)$. 

For $x=(x_{1}, \ldots , x_{N}) \in (\mathbb{C}^{\times})^{N}$ 
we define the reflection operator $T(x) \in {\rm End}(V)$ by 
\begin{align*}
T(x)(v_{a}):=x_{a}^{-1}v_{\overline{a}}, \quad  
T(x)(v_{\overline{a}}):=x_{a}v_{a} \qquad (1 \le a \le N) 
\end{align*}
and set 
\begin{align}
K(\lambda \, | \, x, \beta):=\frac{\lambda\,T(x)+\beta}{\lambda+\beta}.  
\label{eq:def-K}
\end{align}
Then the operator $K(\lambda\,|\,x,\beta)$ satisfies the boundary Yang-Baxter equation: 
\begin{align}
& \label{eq:bYBE}
R_{12}(\lambda_{1}-\lambda_{2})K_{1}(\lambda_{1}\,|\, x, \beta)
R_{21}(\lambda_{1}+\lambda_{2})K_{2}(\lambda_{2}\,|\, x, \beta) \\ 
& \quad{}=
K_{2}(\lambda_{2}\,|\, x, \beta)R_{21}(\lambda_{1}+\lambda_{2})
K_{1}(\lambda_{1}\,|\, x, \beta)R_{12}(\lambda_{1}-\lambda_{2}) 
\nonumber
\end{align}
on $V^{\otimes 2}$, where $K_{j}(\lambda\,|\,x,\beta)$ is the linear operator acting on 
the $j$-th component of $V^{\otimes 2}$ as $K(\lambda\,|\,x,\beta)$. 

For $x=(x_{1}, \ldots , x_{N}) \in (\mathbb{C}^{\times})^{N}$ and 
$y=(y_{1}, \ldots , y_{n}) \in \mathbb{C}^{n}$, 
we define the linear operator $Q_{m}(x\,|\,y) \, (1 \le m \le n)$  acting on $V^{\otimes n}$ by 
\begin{align*}
Q_{m}(x \,|\, y):=&\,
R_{m, m-1}(y_{m}-y_{m-1}-c) \cdots R_{m, 1}(y_{m}-y_{1}-c) 
K_{m}(y_{m}-c/2 \, | \, x, \beta) \\ 
&\times 
R_{1,m}(y_{1}+y_{m}) \cdots R_{m-1,m}(y_{m-1}+y_{m}) \\ 
&\times R_{m,m+1}(y_{m}+y_{m+1}) \cdots R_{m, n}(y_{m}+y_{n}) \\ 
&\times 
K_{m}(y_{m} \, | \, \underbar{$1$}, \alpha)\, 
R_{m, n}(y_{m}-y_{n}) \cdots R_{m, m+1}(y_{m}-y_{m+1}),
\end{align*}
where $\underbar{$1$}:=(1, \ldots , 1)$ and the lower indices of $R$ and $K$ in the right hand side 
signify the position of the components in $V^{\otimes n}$ on which they act. 

\begin{prop}\label{prop:consistency-qkz}
For $1 \le l, m \le n$ we have 
\begin{align*}
Q_{m}(x\,|\,\ldots , y_{l}-c, \ldots)Q_{l}(x \, | \, y)=
Q_{l}(x\,|\,\ldots , y_{m}-c, \ldots)Q_{m}(x \, | \, y). 
\end{align*} 
\end{prop}

\begin{proof}
It follows from \eqref{eq:YBE} and \eqref{eq:bYBE}. 
\end{proof}

Let $f(x\,|\,y)$ be a function on $(\mathbb{C}^{\otimes N})\times\mathbb{C}^{n}$ 
taking values in $V^{\otimes n}$. 
We denote by $\Delta_{m} \, (1 \le m \le n)$ the shift operator with respect to $y_{m}$: 
\begin{align*}
\Delta_{m} f(x \,|\, y):=
f(x \, | \, y_{1}, \ldots , y_{m}-c, \ldots , y_{n}). 
\end{align*}
{}From Proposition \ref{prop:consistency-qkz} 
the following system of difference equations is consistent:  
\begin{align}
\Delta_{m} f(x \,|\, y)=Q_{m}(x \,|\, y)
f(x \,|\, y)  \qquad (m=1, \ldots ,n). 
\label{eq:boundary-rational-qKZ}
\end{align}

\begin{defn}
We call the system \eqref{eq:boundary-rational-qKZ} of difference equations 
the {\it boundary rational qKZ equation}. 
\end{defn}


\section{Compatible differential equations}\label{sec:diff-eqn}

\subsection{Commuting differential operators}

We denote by $e_{ab} \in {\rm End}(V) \, (a, b \in \{1, \ldots , N, \bar{1}, \ldots , \bar{N}\})$ 
the matrix unit acting by $e_{ab}v_{p}=\delta_{bp}v_{a}$. 
In this section, 
for $u \in {\rm End}(V)$ and $1 \le j \le n$, 
we denote by $u^{(j)} \in {\rm End}(V^{\otimes n})$ 
the linear operator acting on the $j$-th component of $V^{\otimes n}$ as $u$. 
For $1 \le a, b \le N$ we set 
\begin{align*}
E_{ab}:=e_{ab}+e_{\overline{a}\overline{b}}, \quad 
\overline{E}_{ab}:=e_{a\overline{b}}+e_{\overline{a}b},   
\end{align*}
and define ${\bf X}_{ab}, {\bf Y}_{ab}, {\bf Z}_{ab} \in {\rm End}(V^{\otimes n})$ by 
\begin{align*}
{\bf X}_{ab}&:=\sum_{1 \le i<j \le n}\left( 
e_{ab}^{(i)}E_{ba}^{(j)}+e_{\overline{b}\overline{a}}^{(i)}E_{ab}^{(j)}
\right), \\ 
{\bf Y}_{ab}&:=\sum_{1 \le i<j \le n}\left( 
e_{a\overline{b}}^{(i)}\overline{E}_{ba}^{(j)}+e_{b\overline{a}}^{(i)}\overline{E}_{ab}^{(j)}
\right), \\ 
{\bf Z}_{ab}&:=\sum_{1 \le i<j \le n}\left( 
e_{\overline{a}b}^{(i)}\overline{E}_{ba}^{(j)}+e_{\overline{b}a}^{(i)}\overline{E}_{ab}^{(j)}
\right).  
\end{align*}
Note that ${\bf Y}_{ab}={\bf Y}_{ba}$ and ${\bf Z}_{ab}={\bf Z}_{ba}$. 

Define the linear operators $A_{a}(y)$ and $B_{a}(x) \, (1 \le a \le N)$ on $V^{\otimes n}$ by
\begin{align*}
A_{a}(y):=\sum_{j=1}^{n}y_{j}(e_{aa}^{(j)}-e_{\overline{a}\overline{a}}^{(j)})+
2\alpha\sum_{j=1}^{n}e_{a\overline{a}}^{(j)}+k\left(-\sum_{p=1}^{a-1}{\bf X}_{pa}+
\sum_{p=a+1}^{N}{\bf X}_{ap}+\sum_{p=1}^{N}{\bf Y}_{ap}
\right)
\end{align*}
and 
\begin{align*}
B_{a}(x)&:=
2\,\frac{\alpha+\beta x_{a}}{x_{a}^{2}-1}\sum_{j=1}^{n}\overline{E}_{aa}^{(j)}+k\,\biggl\{
\sum_{p=1}^{a-1}\frac{x_{a}}{x_{a}-x_{p}}({\bf X}_{ap}+{\bf X}_{pa}) \\ 
&+\sum_{p=a+1}^{N}\frac{x_{p}}{x_{a}-x_{p}}({\bf X}_{ap}+{\bf X}_{pa})+
\sum_{p=1}^{N}\frac{1}{x_{a}x_{p}-1}({\bf Y}_{ap}+{\bf Z}_{ap})
\biggr\}. 
\end{align*}
Set 
\begin{align*}
L_{a}(x \, | \, y):=A_{a}(y)+B_{a}(x) \qquad (1 \le a \le N).  
\end{align*} 
{}By direct calculation we can check the commutativity: 
\begin{lem}\label{lem:hamiltonian}
For $1\le a, b \le N$ we have
\begin{align*}
[A_{a}(y), A_{b}(y)]=0, \qquad 
[L_{a}(x \, | \, y), L_{b}(x \, | \, y)]=0. 
\end{align*} 
\end{lem}

Now define the differential operators $D_{a}(x\,|\,y) \, (1 \le a \le N)$ by 
\begin{align*}
D_{a}(x \, | \, y):=
cx_{a}\frac{\partial}{\partial x_{a}}+L_{a}(x\,|\,y).  
\end{align*}

\begin{prop}
The differential operators $D_{a}(x \, | \, y) \, (a=1, \ldots , N)$ 
commute with each other. 
\end{prop}

\begin{proof}
It follows from Lemma \ref{lem:hamiltonian} and the equality 
$x_{a}\frac{\partial L_{b}}{\partial x_{a}}=x_{b}\frac{\partial L_{a}}{\partial x_{b}}$ 
for $1 \le a, b \le N$.   
\end{proof}

\subsection{Compatibility}

In this subsection  we prove the main theorem: 
\begin{thm}\label{thm:main}
For $1 \le a \le N$ and $1 \le m \le n$ we have 
\begin{align}
[D_{a}(x \,|\, y), \Delta_{m}^{-1}Q_{m}(x \,|\, y)]=0. 
\label{eq:main-commutativity}
\end{align}
Hence the system of equations 
\begin{align}
\left\{ 
\begin{array}{ll}
\Delta_{m}f(x\,|\,y)=Q_{m}(x\,|\,y)f(x\,|\,y) \quad & (1 \le m \le n), \\  & \\
D_{a}(x\,|\,y)f(x\,|\,y)=0 & (1\le a \le N) 
\end{array}
\right.
\label{eq:main}
\end{align} 
is compatible. 
\end{thm}

To prove Theorem \ref{thm:main} we rewrite the linear operators 
$L_{a}(x\,|\,y) \, (1 \le a \le N)$ as follows. 
For $\lambda \in \mathbb{C}^{\times}$ and $\gamma \in \mathbb{C}$, 
we define $I_{a}(\lambda \, | \, \gamma) \in {\rm End}(V)$ by 
\begin{align*}
I_{a}(\lambda \, | \, \gamma):=
\gamma(e_{aa}-e_{\overline{a}\overline{a}})+2\frac{\alpha+\beta \lambda}{\lambda^2-1}e_{\overline{a}a}+
2\frac{\alpha+\beta \lambda^{-1}}{1-\lambda^{-2}}e_{a\overline{a}}
\quad (1 \le a \le N). 
\end{align*}
For $x=(x_{1}, \ldots , x_{N}) \in (\mathbb{C}^{\times})^{N}$ 
we define $M_{a}(x) \in {\rm End}(V^{\otimes 2}) \, (1 \le a \le N)$ by
\begin{align*}
M_{a}(x)&:= 
\frac{2k}{x_{a}-x_{a}^{-1}}
(x_{a}e_{a\overline{a}}+x_{a}^{-1}e_{\overline{a}a})\otimes 
(e_{a\overline{a}}+e_{\overline{a}a}) \\ 
&+k\sum_{\begin{subarray}{c}p=1 \\ p\not=a\end{subarray}}^{N}\left( 
\frac{x_{a}}{x_{a}-x_{p}}U_{ap}+\frac{x_{p}}{x_{a}-x_{p}}U_{pa}+
\frac{x_{a}x_{p}}{x_{a}x_{p}-1}J_{ap}+
\frac{1}{x_{a}x_{p}-1}K_{ap}
\right), 
\end{align*}
where
\begin{align*}
U_{ab}&:=e_{ab}\otimes E_{ba}+e_{\overline{b}\overline{a}} \otimes E_{ab}, \qquad 
J_{ab}:=e_{a\overline{b}}\otimes\overline{E}_{ba}+e_{b\overline{a}}\otimes\overline{E}_{ab}, \\
K_{ab}&:=e_{\overline{a}b}\otimes\overline{E}_{ba}+e_{\overline{b}a}\otimes\overline{E}_{ab}. 
\end{align*}
Note that $J_{ab}=J_{ba}$ and $K_{ab}=K_{ba}$. 
Then we have 
\begin{align*}
L_{a}(x \, | \, y)=\sum_{j=1}^{n}I_{a}(x_{a} \, | \, y_{j})^{(j)}+
\sum_{1 \le i<j \le n}M_{a}(x)^{(i,j)} 
\end{align*}
for $1 \le a \le N$, where $M_{a}(x)^{(i,j)}$ is the linear operator 
acting on the tensor product of the $i$-th and the $j$-th component of $V^{\otimes n}$ as $M_{a}(x)$.

\begin{lem}\label{lem:comm-IM}
For $x=(x_{1}, \ldots , x_{N})$ and $1 \le a \le N$ we have 
\begin{align*}
& 
R_{12}(y_{1}-y_{2}) \left( 
I_{a}(x_{a} | y_{1})^{(1)}+I_{a}(x_{a} | y_{2})^{(2)}+M_{a}(x)^{(1,2)} 
\right) 
R_{12}(y_{1}-y_{2})^{-1} \\ 
&=
I_{a}(x_{a} | y_{1})^{(1)}+I_{a}(x_{a} | y_{2})^{(2)}+M_{a}(x)^{(2,1)} 
\end{align*} 
on $V^{\otimes 2}$, where $M_{a}(x)^{(2,1)}=PM_{a}(x)P$. 
\end{lem}

\begin{proof}
Note that if $h \in {\rm End}(V^{\otimes 2})$ is symmetric, i.e. $PhP=h$, 
then $h$ commutes with the $R$-matrix. 
We extract symmetric parts from the operator 
$I_{a}(x_{a} | y_{1})^{(1)}+I_{a}(x_{a} | y_{2})^{(2)}+M_{a}(x)^{(1,2)}$ 
as follows. 
In the following we enclose symmetric parts in a square bracket $[\quad]$. 
First we have 
\begin{align*}
& 
I_{a}(x_{a} | y_{1})^{(1)}+I_{a}(x_{a} | y_{2})^{(2)} \\
&=\left[ 
y_{1}\sum_{j=1}^{2}(e_{aa}^{(j)}-e_{\bar{a}\bar{a}}^{(j)})+
2\frac{\alpha+\beta x_{a}}{x_{a}^{2}-1}\sum_{j=1}^{2}e_{\bar{a}a}^{(j)}+
2\frac{\alpha+\beta x_{a}^{-1}}{1-x_{a}^{-2}}\sum_{j=1}^{2}e_{a\bar{a}}^{(j)}
\right] \\ 
&+(y_{2}-y_{1})(e_{aa}^{(2)}-e_{\bar{a}\bar{a}}^{(2)}). 
\end{align*}
Note that $U_{ap}+U_{pa}$ and $J_{ap}+K_{ap}$ are symmetric. 
Using 
\begin{align*}
U_{ap}&=\left[
e_{\bar{p}\bar{a}}^{(1)}e_{ap}^{(2)}+e_{ap}^{(1)}e_{\bar{p}\bar{a}}^{(2)}
\right]+
e_{\bar{p}\bar{a}}^{(1)}e_{\bar{a}\bar{p}}^{(2)}+e_{ap}^{(1)}e_{pa}^{(2)}, \\ 
J_{ap}&=\left[ 
e_{p\bar{a}}^{(1)}e_{a\bar{p}}^{(2)}+e_{a\bar{p}}^{(1)}e_{p\bar{a}}^{(2)}
\right]+
e_{p\bar{a}}^{(1)}e_{\bar{a}p}^{(2)}+e_{a\bar{p}}^{(1)}e_{\bar{p}a}^{(2)},  
\end{align*}
we divide $M_{a}(x)^{(1,2)}$ as 
\begin{align*}
& 
M_{a}(x)^{(1,2)}=\frac{2k}{x_{a}-x_{a}^{-1}}
\left[
x_{a}e_{a\bar{a}}^{(1)}e_{a\bar{a}}^{(2)}+x_{a}^{-1}e_{\bar{a}a}^{(1)}e_{\bar{a}a}^{(2)}+
x_{a}^{-1}(e_{a\bar{a}}^{(1)}e_{\bar{a}a}^{(2)}+e_{\bar{a}a}^{(1)}e_{a\bar{a}}^{(2)})
\right] \\ 
&+k\sum_{\begin{subarray}{c}p=1 \\ p\not=a\end{subarray}}^{N}\left[
\frac{x_{p}}{x_{a}-x_{p}}(U_{ap}+U_{pa})+\frac{1}{x_{a}x_{p}-1}(J_{ap}+K_{ap})
\right] \\ 
&+k\left[\sum_{\begin{subarray}{c}p=1 \\ p\not=a\end{subarray}}^{N}(
e_{\bar{p}\bar{a}}^{(1)}e_{ap}^{(2)}+e_{ap}^{(1)}e_{\bar{p}\bar{a}}^{(2)}+
e_{p\bar{a}}^{(1)}e_{a\bar{p}}^{(2)}+e_{a\bar{p}}^{(1)}e_{p\bar{a}}^{(2)})-
(e_{aa}^{(1)}e_{aa}^{(2)}+e_{\bar{a}\bar{a}}^{(1)}e_{\bar{a}\bar{a}}^{(2)})
\right] \\ 
&+k\sideset{}{'}\sum_{p}(e_{ap}^{(1)}e_{pa}^{(2)}+e_{p\bar{a}}^{(1)}e_{\bar{a}p}^{(2)}), 
\end{align*}
where $\textstyle \sum_{p}'$ is a sum over all indices $p \in \{1, \ldots , N, \bar{1}, \ldots , \bar{N}\}$. 
Thus we find 
\begin{align}
& \label{eq:split-symm}
I_{a}(x_{a} | y_{1})^{(1)}+I_{a}(x_{a} | y_{2})^{(2)}+M_{a}(x)^{(1,2)} \\ 
&=\left[ (\hbox{symmetric part}) \right]+
(y_{2}-y_{1})(e_{aa}^{(2)}-e_{\bar{a}\bar{a}}^{(2)})+
k\sideset{}{'}\sum_{p}(e_{ap}^{(1)}e_{pa}^{(2)}+e_{p\bar{a}}^{(1)}e_{\bar{a}p}^{(2)}). 
\nonumber
\end{align}
Now we make use of the intertwining property of the $R$-matrix: 
\begin{align*}
R(\lambda)\left( \lambda \cdot 1 \otimes e_{lm}+k\sideset{}{'}\sum_{p}e_{pm}\otimes e_{lp} \right)&=
\left( \lambda \cdot 1 \otimes e_{lm}+k\sideset{}{'}\sum_{p}e_{lp}\otimes e_{pm} \right)R(\lambda), \\ 
R(\lambda)\left( \lambda \cdot 1 \otimes e_{lm}-k\sideset{}{'}\sum_{p}e_{lp}\otimes e_{pm} \right)&=
\left( \lambda \cdot 1 \otimes e_{lm}-k\sideset{}{'}\sum_{p}e_{pm}\otimes e_{lp} \right)R(\lambda) 
\end{align*}
for any $l, m \in \{1, \ldots , N, \bar{1}, \ldots , \bar{N}\}$. 
For $g \in {\rm GL}(V^{\otimes 2})$ we denote by $\mathop{\rm Ad}(g)$ 
the adjoint action $h \mapsto ghg^{-1}$ on ${\rm End}(V^{\otimes 2})$.
Then we see that 
\begin{align*}
&
\mathop{\rm Ad}(R_{12}(y_{1}-y_{2}))
\left( (y_{2}-y_{1})(e_{aa}^{(2)}-e_{\bar{a}\bar{a}}^{(2)})+
k\sideset{}{'}\sum_{p}(e_{ap}^{(1)}e_{pa}^{(2)}+e_{p\bar{a}}^{(1)}e_{\bar{a}p}^{(2)}) 
\right) \\
&=(y_{2}-y_{1})(e_{aa}^{(2)}-e_{\bar{a}\bar{a}}^{(2)})+
k\sideset{}{'}\sum_{p}(e_{pa}^{(1)}e_{ap}^{(2)}+e_{\bar{a}p}^{(1)}e_{p\bar{a}}^{(2)}) \\ 
&=
(y_{2}-y_{1})(e_{aa}^{(2)}-e_{\bar{a}\bar{a}}^{(2)})+2k\,e_{\bar{a}a}^{(1)}e_{a\bar{a}}^{(2)} \\ 
&+k\biggl\{\sum_{\begin{subarray}{c}p=1 \\ p\not=a\end{subarray}}^{N}(
e_{pa}^{(1)}e_{ap}^{(2)}+e_{\bar{p}a}^{(1)}e_{a\bar{p}}^{(2)}+
e_{\bar{a}p}^{(1)}e_{p\bar{a}}^{(2)}+e_{\bar{a}\bar{p}}^{(1)}e_{\bar{p}\bar{a}}^{(2)}
)+(e_{aa}^{(1)}e_{aa}^{(2)}+e_{\bar{a}\bar{a}}^{(1)}e_{\bar{a}\bar{a}}^{(2)})
\biggr\}. 
\end{align*}
Adding this to the symmetric part in \eqref{eq:split-symm} 
we obtain $I_{a}(x_{a} | y_{1})^{(1)}+I_{a}(x_{a} | y_{2})^{(2)}+M_{a}(x)^{(2,1)}$. 
\end{proof}

\begin{proof}[Proof of Theorem \ref{thm:main}] 
We split the operator $Q_{m}(x\,|\,y)$ into three parts: 
\begin{align}
Q_{m}(x\,|\,y)=Q_{m}'(y)\,K_{m}(y_{m}-c/2\,|\,x, \beta)\,Q_{m}''(y),  
\label{eq:Q-split}
\end{align}
where 
\begin{align*}
Q_{m}'(y):=R_{m, m-1}(y_{m}-y_{m-1}-c) \cdots R_{m, 1}(y_{m}-y_{1}-c)  
\end{align*}
and $Q_{m}''(y)$ is determined by \eqref{eq:Q-split}. 
Then the equality \eqref{eq:main-commutativity} is equivalent to 
\begin{align}
& \label{eq:proof-comm}
cx_{a}\left(\frac{\partial}{\partial x_{a}}K_{m}(y_{m}-c/2\,|\,x,\beta)\right)
K_{m}(y_{m}-c/2\,|\,x,\beta)^{-1} \\ 
&+\mathop{\rm Ad}(Q_{m}'(y))(L_{a}(x\,| \ldots , y_{m}-c, \ldots)) 
\nonumber \\ 
&-\mathop{\rm Ad}\left(K_{m}(y_{m}-c/2\,|\,x, \beta)Q_{m}''(y)\right) 
(L_{a}(x\,|\ldots , y_{m}, \ldots))=0. 
\nonumber 
\end{align}

The first term of \eqref{eq:proof-comm} is equal to 
\begin{align}
\frac{c(y_{m}-c/2)}{(y_{m}-c/2)^{2}-\beta^2}\left\{ 
(y_{m}-c/2)(e_{aa}^{(m)}-e_{\bar{a}\bar{a}}^{(m)})-
\beta(x_{a}e_{a\bar{a}}^{(m)}-x_{a}^{-1}e_{\bar{a}a}^{(m)})
\right\}. 
\label{eq:proof-comm-1} 
\end{align}
{}From Lemma \ref{lem:comm-IM} we see that the second term is equal to 
\begin{align}
\sum_{\begin{subarray}{c} j=1 \\ j\not=m \end{subarray}}^{n}I_{a}(x_{a}\,|\,y_{j})^{(j)}+
I_{a}(x_{a}\,|\,y_{m}-c)^{(m)}+
\sum_{\begin{subarray}{c} j=1 \\ j\not=m \end{subarray}}^{n}M_{a}(x)^{(m,j)}+
\sum_{\begin{subarray}{c} 1\le i<j\le n\\ i,j\not=m \end{subarray}}^{n}M_{a}(x)^{(i,j)}. 
\label{eq:proof-comm-2} 
\end{align}
Let us calculate the third term. 
Using Lemma \ref{lem:comm-IM} and 
\begin{align*}
& 
\mathop{\rm Ad}(K(\gamma\,|\,\underbar{$1$}, \alpha))(I_{a}(\lambda\,|\,\gamma))=
I_{a}(\lambda\,|\,-\gamma),  \\ 
& 
\mathop{\rm Ad}(K_{2}(\gamma\,|\,\underbar{$1$}, \alpha))(M_{a}(x)^{(1,2)})=M_{a}(x)^{(1,2)}, 
\end{align*}
we obtain 
\begin{align*}
& 
\mathop{\rm Ad}(Q_{m}''(y))(L_{a}(x\,|\,y)) \\ 
&=\sum_{\begin{subarray}{c} j=1 \\ j\not=m \end{subarray}}^{n}I_{a}(x_{a}\,|\,y_{j})^{(j)}+
I_{a}(x_{a}\,|-y_{m})^{(m)}+
\sum_{\begin{subarray}{c} j=1 \\ j\not=m \end{subarray}}^{n}M_{a}(x)^{(m,j)}+
\sum_{\begin{subarray}{c} 1\le i<j\le n\\ i,j\not=m \end{subarray}}^{n}M_{a}(x)^{(i,j)}. 
\end{align*}
What remains is to calculate the image of the second and the third term above by 
the operator $\mathop{\rm Ad}(K_{m}(y_{m}-c/2\,|\,x,\beta))$. 
By direct calculation we find 
\begin{align*}
& 
\mathop{\rm Ad}(K(\gamma-c/2\,|\,x,\beta))(I_{a}(x_{a}\,|-\gamma)) \\ 
&=I_{a}(x_{a}\,|\,\gamma)+\frac{c\beta}{(\gamma-c/2)^{2}-\beta^{2}}\left\{ 
\beta(e_{aa}-e_{\bar{a}\bar{a}})+(\gamma-c/2)(x_{a}^{-1}e_{\bar{a}a}-x_{a}e_{a\bar{a}})
\right\}
\end{align*}
and 
\begin{align*}
\mathop{\rm Ad}(K_{1}(\gamma\,|\,x,\beta))(M_{a}(x)^{(1,2)})=M_{a}(x)^{(1,2)}. 
\end{align*}
Using these formulas we see that the third term of \eqref{eq:proof-comm} is equal to 
\begin{align}
\label{eq:proof-comm-3}
(-1)\times \biggl( & 
\sum_{j=1}^{n}I_{a}(x_{a}\,|\,y_{j})^{(j)}+
\sum_{\begin{subarray}{c} j=1 \\ j\not=m \end{subarray}}^{n}M_{a}(x)^{(m,j)}+
\sum_{\begin{subarray}{c} 1\le i<j\le n\\ i,j\not=m \end{subarray}}^{n}M_{a}(x)^{(i,j)} \\ 
&+\frac{c\beta}{(y_{m}-c/2)^{2}-\beta^{2}}\left\{ 
\beta(e_{aa}^{(m)}-e_{\bar{a}\bar{a}}^{(m)})+(y_{m}-c/2)(x_{a}^{-1}e_{\bar{a}a}^{(m)}-x_{a}e_{a\bar{a}}^{(m)})
\right\}
\biggr).  
\nonumber
\end{align}
The sum of \eqref{eq:proof-comm-1}, \eqref{eq:proof-comm-2} and \eqref{eq:proof-comm-3} is zero 
and this completes the proof. 
\end{proof}


\section{The bispectral qKZ equation and its degeneration}\label{sec:bispectral}

\subsection{The double affine Hecke algebra of type $(C_{n}^{\vee},C_{n})$}

Here we recall the definition and some properties of the double affine Hecke algebra 
of type $(C_{n}^{\vee},C_{n})$ \cite{S}. 
We denote by 
$\mathbb{F}:=\mathbb{C}(q^{1/2}, t^{1/2}, t_{0}^{1/2}, t_{n}^{1/2}, u_{0}^{1/2}, u_{n}^{1/2})$ 
the coefficient field. 

\begin{defn}
The {\it double affine Hecke algebra} $\DAHA$ of type $(C_{n}^{\vee},C_{n})$ is the unital associative 
$\mathbb{F}$-algebra generated by $X_{i}^{\pm 1} \, (1 \le i \le n)$ and $T_{i} \, (0 \le i \le n)$ 
satisfying the following relations:

(i)\, quadratic Hecke relations
\begin{align*}
(T_{i}-t_{i}^{1/2})(T_{i}+t_{i}^{-1/2})=0 \qquad (0 \le i \le n), 
\end{align*} 
where $t_{i}^{1/2}:=t^{1/2}$ for $1 \le i<n$. 

(ii)\, braid relations
\begin{align*}
& 
T_{i}T_{i+1}T_{i}T_{i+1}=T_{i+1}T_{i}T_{i+1}T_{i} \quad (i=0, n-1), \\
& 
T_{i}T_{i+1}T_{i}=T_{i+1}T_{i}T_{i+1} \quad (1\le i<n), \qquad 
T_{i}T_{j}=T_{j}T_{i} \quad (|i-j|\ge 2). 
\end{align*}

(iii)\, relations between $X$ and $T$
\begin{align*}
& 
X_{i}X_{j}=X_{j}X_{i} \quad (\forall{i, j}), \qquad 
T_{i}X_{j}=X_{j}T_{i} \quad (|i-j| \ge 2 \,\,\hbox{or}\,\,(i,j)=(n,n-1)), \\ 
& 
T_{i}X_{i}T_{i}=X_{i+1} \quad (1 \le i \le n-1), \qquad 
X_{n}^{-1}T_{n}^{-1}=T_{n}X_{n}+(u_{n}^{1/2}-u_{n}^{-1/2}), \\
&
q^{-1/2}T_{0}^{-1}X_{1}=q^{1/2}X_{1}^{-1}T_{0}+(u_{0}^{1/2}-u_{0}^{-1/2}). 
\end{align*}
\end{defn}

Noumi found the polynomial representation of $\DAHA$ given as follows \cite{N}. 
Let $W=\langle s_{0}, \ldots , s_{n} \rangle$ be the affine Weyl group of type $C_{n}$. 
The group $W$ acts on the Laurent polynomial ring 
$\mathbb{F}[X^{\pm 1}]=\mathbb{F}[X_{1}^{\pm 1}, \ldots , X_{n}^{\pm 1}]$ by 
\begin{align}
(s_{0}f)(X)&=f(qX_{1}^{-1}, X_{2}, \ldots , X_{n}), 
\label{eq:weyl-action-1} \\ 
(s_{i}f)(X)&=f(\ldots , X_{i+1}, X_{i}, \ldots) \qquad (1 \le i <n), 
\label{eq:weyl-action-2} \\
(s_{n}f)(X)&=f(X_{1}, \ldots , X_{n-1}, X_{n}^{-1}).
\label{eq:weyl-action-3}
\end{align}
Define the $\mathbb{F}$-linear operators $\widehat{T}_{i} \,\, (0\le i \le n)$ 
on $\mathbb{F}[X^{\pm 1}]$: 
\begin{align*}
\widehat{T}_{0}&:=
t_{0}^{1/2}+t_{0}^{-1/2}
\frac{(1-q^{1/2}t_{0}^{1/2}u_{0}^{1/2}X_{1}^{-1})(1+q^{1/2}t_{0}^{1/2}u_{0}^{-1/2}X_{1}^{-1})}
     {1-qX_{1}^{-2}}
(s_{0}-1), \\ 
\widehat{T}_{i}&:=
t^{1/2}+t^{-1/2}\frac{1-tX_{i}/X_{i+1}}{1-X_{i}/X_{i+1}}(s_{i}-1) \qquad (1 \le i <n), \\ 
\widehat{T}_{n}&:=
t_{n}^{1/2}+t_{n}^{-1/2}
\frac{(1-t_{n}^{1/2}u_{n}^{1/2}X_{n})(1+t_{n}^{1/2}u_{n}^{-1/2}X_{n})}{1-X_{n}^{2}}(s_{n}-1).
\end{align*}
Then the map $T_{i} \mapsto \widehat{T}_{i}$ and $X_{i} \mapsto X_{i}$ (left multiplication) 
gives a representation of $\DAHA$ on $\mathbb{F}[X^{\pm 1}]$. 
It is faithful and hence $\DAHA$ is isomorphic to the $\mathbb{F}$-subalgebra 
of the smashed product algebra $\mathbb{F}(X) \# W$ 
generated by the difference operators $\widehat{T}_{i} \, (0 \le i \le n)$ and 
the multiplication operators $f(\underbar{X}) \in \mathbb{F}[X^{\pm 1}]$ 
(see, e.g., Section 2.1 in \cite{MS} for the definition of the smashed product algebra). 
Hereafter we identify $\DAHA$ with the subalgebra of $\mathbb{F}(X) \# W$. 

The subalgebra $H_{0}$ generated by $T_{i} \, (1 \le i \le n)$ is isomorphic to 
the Hecke algebra of type $C_{n}$. 
Denote by $W_{0}:=\langle s_{1}, \ldots , s_{n} \rangle$ the finite Weyl group of type $C_{n}$. 
Let $w=s_{j_{1}} \cdots s_{j_{r}}$ be a reduced expression of $w \in W_{0}$.
Then the element $T_{w}:=T_{j_{1}} \cdots T_{j_{r}}$ is well-defined for $w \in W_{0}$.  
The set $\{ T_{w}\}_{w \in W_{0}}$ gives a basis of $H_{0}$.  

Set 
\begin{align*}
Y_{i}:=T_{i}\cdots T_{n-1}(T_{n} \cdots T_{0}) T_{1}^{-1} \cdots T_{i-1}^{-1} \quad (1 \le i \le n). 
\end{align*}
They satisfy 
\begin{align*}
& 
Y_{i}Y_{j}=Y_{j}Y_{i} \quad (\forall{i, j}), \qquad 
T_{i}Y_{j}=Y_{j}T_{i} \quad (|i-j| \ge 2 \,\,\hbox{or}\,\,(i,j)=(n,n-1)), \\  
& 
T_{i}Y_{i+1}T_{i}=Y_{i} \quad (1 \le i \le n-1), \qquad 
T_{n}^{-1}Y_{n}=Y_{n}^{-1}T_{n}+(t_{0}^{1/2}-t_{0}^{-1/2})
\end{align*} 
and 
\begin{align*}
q^{-1/2}Y_{1}^{-1}U_{n}^{-1}=q^{1/2}U_{n}Y_{1}+(u_{0}^{1/2}-u_{0}^{-1/2}), 
\end{align*}
where $U_{n}:=X_{1}^{-1}T_{0}Y_{1}^{-1}$. 
The subalgebra $H$ generated by $T_{i} \, (1 \le i \le n)$ and $Y_{i}^{\pm 1} \, (1 \le i \le n)$ is 
called the {\it affine Hecke algebra} of type $C_{n}$. 

Let $* : \mathbb{F} \to \mathbb{F}$ be the $\mathbb{C}$-algebra involution 
defined by $(t_{0}^{1/2})^{*}=u_{n}^{1/2}$ and 
the other parameters $q^{1/2}, t^{1/2}, t_{n}^{1/2}, u_{0}^{1/2}$ are fixed.  
It uniquely extends to the $\mathbb{C}$-algebra anti-involution on $\DAHA$ such that 
\begin{align*}
T_{0}^{*}=U_{n}, \quad 
T_{i}^{*}=T_{i}, \quad 
X_{i}^{*}=Y_{i}^{-1}, \quad 
Y_{i}^{*}=X_{i}^{-1}  \qquad (1 \le i \le n). 
\end{align*}
The anti-involution $*$ is called duality anti-involution. 
For a Laurent polynomial $f$ of $n$-variables with the coefficients in $\mathbb{F}$, 
we define $f^{\diamond}$ by the equality $f^{\diamond}(Y)=(f(X))^{*}$. 

Define the elements $\widetilde{S}_{i} \, (0 \le i \le n)$ in $\mathbb{F}(X) \# W$ by 
\begin{align*}
\widetilde{S}_{0}&:=t_{0}^{-1/2}
(1-q^{1/2}t_{0}^{1/2}u_{0}^{1/2}X_{1}^{-1})(1+q^{1/2}t_{0}^{1/2}u_{0}^{-1/2}X_{1}^{-1})s_{0}, \\ 
\widetilde{S}_{i}&:=t^{-1/2}(1-tX_{i}/X_{i+1})s_{i} \quad (1 \le i<n), \\ 
\widetilde{S}_{n}&:=t_{n}^{-1/2}
(1-t_{n}^{1/2}u_{n}^{1/2}X_{n})(1+t_{n}^{1/2}u_{n}^{-1/2}X_{n})s_{n}. 
\end{align*}
In fact they belong to the subalgebra $\DAHA$ since we have 
\begin{align}
\widetilde{S}_{0}&=(1-qX_{1}^{-2})T_{0}-(t_{0}^{1/2}-t_{0}^{-1/2})-(u_{0}^{1/2}-u_{0}^{-1/2})q^{1/2}X_{1}^{-1}, 
\label{eq:intertwiner-generator} \\ 
\widetilde{S}_{i}&=(1-X_{i}/X_{i+1})T_{i}-(t^{1/2}-t^{-1/2}) \quad (1 \le i <n), 
 \nonumber \\ 
\widetilde{S}_{n}&=(1-X_{n}^{2})T_{n}-(t_{n}^{1/2}-t_{n}^{-1/2})-(u_{n}^{1/2}-u_{n}^{-1/2})X_{n}.   
\nonumber 
\end{align}
The elements $\widetilde{S}_{i}$ and 
their dual $\widetilde{S}_{i}^{*} \, (0 \le i \le n)$ will play a fundamental role 
in the construction of the bispectral qKZ equation.

\subsection{The bispectral qKZ equation}

Here we construct the bispectral qKZ equation of type $(C_{n}^{\vee}, C_{n})$. 
See \cite{MS} for the details in the case of type $GL_{n}$. 

Hereafter we set the parameters $q^{1/2}, t^{1/2}, \ldots$ to generic complex values 
and consider $\mathbb{H}$ as a $\mathbb{C}$-algebra. 
Then we have the Poincar\'e-Birkhoff-Witt (PBW) decomposition of $H_{0}$ and $\DAHA$ 
as $\mathbb{C}$-vector spaces: 
\begin{align*}
H \simeq H_{0} \otimes \mathbb{C}[Y^{\pm 1}], \qquad 
\DAHA \simeq 
\mathbb{C}[X^{\pm 1}] \otimes H_{0} \otimes \mathbb{C}[Y^{\pm 1}]. 
\end{align*}

Set $\mathbb{L}:=\mathbb{C}[X^{\pm 1}] \otimes \mathbb{C}[Y^{\pm 1}]$. 
The DAHA $\DAHA$ has $\mathbb{L}$-module structure defined by 
\begin{align}
(f \otimes g).h=f(x) \, h \, g(y) \qquad 
(f\otimes g \in \mathbb{L}, \, h \in \DAHA). 
\label{eq:L-mod-str}
\end{align}
We consider $\mathbb{L}$ as the ring of regular functions on $T \times T$, 
where $T$ is the $n$-dimensional torus $T:=(\mathbb{C}^{\times})^{n}$. 
{}From the PBW decomposition, any element of $\DAHA$ can be regarded as an $H_{0}$-valued 
regular function on $T \times T$. 
Let $\mathbb{K}$ be the field of meromorphic functions on $T \times T$. 
Then $H_{0}^{\mathbb{K}}:=\mathbb{K} \otimes_{\mathbb{L}} \DAHA$ is a left $\mathbb{K}$-module 
of $H_{0}$-valued meromorphic functions on $T \times T$. 
Any element $F \in H_{0}^{\mathbb{K}}$ is uniquely written in the form 
$F=\sum_{w \in W_{0}}f_{w}.T_{w} \,\, (f_{w} \in \mathbb{K})$. 

Denote the translations in $W$ by 
\begin{align*}
\epsilon_{i}:=s_{i} \cdots s_{n-1} (s_{n} \cdots s_{0}) s_{1} \cdots s_{i-1} \qquad (1 \le i \le n).  
\end{align*}
Then $W$ is a semi-direct product $W \simeq W_{0} \ltimes \Gamma$ of  
the finite Weyl group $W_{0}:=\langle s_{1}, \ldots , s_{n} \rangle$ and 
the lattice $\Gamma \simeq \mathbb{Z}^{n}$ generated by $\epsilon_{i} \, (1 \le i \le n)$. 
Define the involution ${}^{\diamond} \, : \, W \to W$ by 
$w_{0}^{\diamond}=w_{0}$ for $w_{0} \in W_{0}$ and 
$\epsilon_{i}^{\diamond}=\epsilon_{i}^{-1} \, (1 \le i \le n)$. 
Then $W\times W$ acts on $\mathbb{L}$ by 
$(w, w')(f(X)\otimes g(Y))=(wf)(X)\otimes (w'^{\diamond}g)(Y)$. 
The action naturally extends to that on $\mathbb{K}$. 
Now we define the action of $W \times W$ on $H_{0}^{\mathbb{K}}$ 
by $(w, w')(F):=\sum_{w \in W_{0}}(w, w')(f_{w}).T_{w}$ for 
$F=\sum_{w \in W_{0}}f_{w}.T_{w} \in H_{0}^{\mathbb{K}}$. 

Let $w=s_{j_{1}} \cdots s_{j_{l}}$ be a reduced expression of $w \in W$. 
Then the element $\widetilde{S}_{w}:=\widetilde{S}_{j_{1}} \cdots \widetilde{S}_{j_{l}} \in \DAHA$ is well-defined.  
We denote by $d_{w}(X)$ the Laurent polynomial uniquely determined by 
the equality $\widetilde{S}_{w}=d_{w}(X)w$ in $\mathbb{F}(X)\# W$.  

For $(w, w') \in W \times W$, consider the $\mathbb{C}$-linear endomorphism 
$\widetilde{\sigma}_{(w, w')}$ on $\DAHA$ defined by 
\begin{align*}
\widetilde{\sigma}_{(w, w')}(h):=\widetilde{S}_{w}\,h\,\widetilde{S}_{w'}^{*}.   
\end{align*}
Then we have 
\begin{align}
\widetilde{\sigma}_{(w, w')}(f.h)=(w, w')(f).\widetilde{\sigma}_{(w, w')}(h) \qquad 
(f \in \mathbb{L}, \, h \in \DAHA).  
\label{eq:sigma-tilde}
\end{align}
The map $\widetilde{\sigma}_{(w, w')}$ extends to the $\mathbb{C}$-linear endomorphism on $H_{0}^{\mathbb{K}}$ 
satisfying \eqref{eq:sigma-tilde} for all $f \in \mathbb{K}$ and $h \in \DAHA$. 

We define $\tau(w, w') \in {\rm End}_{\mathbb{C}}(H_{0}^{\mathbb{K}}) \, (w, w'\in W)$ by 
\begin{align*}
\tau(w, w')(F):=d_{w}(X)^{-1}d_{w'}^{\diamond}(Y)^{-1}. 
\widetilde{\sigma}_{(w, w')}(F) \qquad (F \in H_{0}^{\mathbb{K}}).   
\end{align*}
Using the equality 
\begin{align*}
d_{w_{1}}(X)^{-1}(w_{1}d_{w_{2}})(X)^{-1}\widetilde{S}_{w_{1}}\widetilde{S}_{w_{2}}=
w_{1}w_{2}=d_{w_{1}w_{2}}(X)^{-1}\widetilde{S}_{w_{1}w_{2}}, 
\end{align*}
we see that $\tau$ is a group homomorphism $\tau: W\times W \to {\rm GL}_{\mathbb{C}}(H_{0}^{\mathbb{K}})$. 
{}From the definition of $\tau$, the operator 
\begin{align*}
C(w, w'):=\tau(w, w') \cdot (w, w')^{-1}  
\end{align*}
acting on $H_{0}^{\mathbb{K}}$ is $\mathbb{K}$-linear. 

Now consider the equation 
\begin{align*}
\tau(\mu, \nu)(F)=F \qquad (\forall{\mu, \nu} \in \Gamma) 
\end{align*}
for $F \in H_{0}^{\mathbb{K}}$. 
It is rewritten as 
\begin{align}
C(\mu, \nu)F(\mu X \, |\, \nu^{-1}Y)=F(X \,|\, Y) \qquad (\forall{\mu, \nu} \in \Gamma), 
\label{eq:bispectral-qKZ} 
\end{align}
where $\mu Z:=(q^{\mu_{1}}Z_{1}, \ldots , q^{\mu_{n}}Z_{n})$ for 
$Z=(Z_{1}, \ldots , Z_{n})$ and $\mu=\epsilon_{1}^{\mu_{1}} \cdots \epsilon_{n}^{\mu_{n}} \in \Gamma$. 
Thus \eqref{eq:bispectral-qKZ} is a system of linear $q$-difference equations. 
Since $\tau$ is a group homomorphism, the system is holonomic. 
 
\begin{definition}
We call the holonomic system of $q$-difference equations $\eqref{eq:bispectral-qKZ}$ 
the {\it bispectral qKZ equation of type} $(C_{n}^{\vee}, C_{n})$.   
\end{definition}
 
The bispectral qKZ equation \eqref{eq:bispectral-qKZ} essentially consists of 
the two systems:
\begin{align}
& 
C(\epsilon_{a}, 1)F(\ldots , qX_{a}, \ldots |Y)=
F(\ldots , X_{a}, \ldots |\, Y) \quad (1 \le a \le n), 
\label{eq:QAKZ-X} \\
& 
C(1, \epsilon_{m})F(X|\ldots , q^{-1}Y_{m}, \ldots )=
F(X\,|\ldots , Y_{m}, \ldots) \quad (1 \le m \le n).
\label{eq:QAKZ-Y} 
\end{align} 
The system \eqref{eq:QAKZ-X} is called the 
{\it quantum affine Knizhnik-Zamolodchikov} (QAKZ) {\it equation of type} $C_{n}$. 
See Section 1.3.6 of \cite{C2} for construction of the QAKZ equation 
associated with arbitrary root system. 
The system \eqref{eq:QAKZ-Y} is dual of \eqref{eq:QAKZ-X}. 
In the rest of this subsection  we compute the operator $C(1, \epsilon_{m}) \, (1 \le m \le n)$ explicitly.  
 
Denote by $H^{*}$ the subalgebra of $\mathbb{H}$ generated by 
$X_{i}^{\pm 1} \, (1 \le i \le n)$ and $T_{i} \, (1 \le i \le n)$. 
The duality anti-involution $*$ gives the isomorphism $H \simeq H^{*}$. 
We define the anti-algebra homomorphism 
$\eta_{R}: H^{*} \to {\rm End}_{\mathbb{K}}(H_{0}^{\mathbb{K}})$ 
by 
\begin{align*}
\eta_{R}(A)\left(\sum_{w\in W_{0}}f_{w}.T_{w}\right):=\sum_{w \in W_{0}}f_{w}.(T_{w}A),   
\end{align*}
where $A \in H^{*}$ and $f_{w} \in \mathbb{K} \, (w \in W_{0})$. 
Applying the dual anti-involution to \eqref{eq:intertwiner-generator}, 
we obtain explicit formulas for $\widetilde{S}^{*}_{i}$ in terms of 
$T_{i} \, (1 \le i \le n), U_{n}$ and $Y_{i}^{\pm 1} \, (1 \le i \le n)$. 
Then we see that $C(1, s_{i}) \, (0 \le i \le n)$ are given by 
\begin{align*}
C(1, s_{i})=\left\{ 
\begin{array}{ll}
\mathcal{K}_{0}(Y_{1}) & (i=0), \\
\mathcal{R}_{i}(Y_{i+1}/Y_{i}) & (1\le i<n), \\
\mathcal{K}_{n}(Y_{n}) & (i=n), 
\end{array}
\right.
\end{align*}
where $\mathcal{K}_{0}(Y), \mathcal{R}_{i}(Y) \, (1 \le i<n)$ and $\mathcal{K}_{n}(Y)$ 
are defined by 
\begin{align*}
\mathcal{K}_{0}(Y)&:=\frac{u_{n}^{1/2}}{(1-q^{1/2}u_{0}^{1/2}u_{n}^{1/2}Y)(1+q^{1/2}u_{0}^{-1/2}u_{n}^{1/2}Y)} \\ 
&\times\left\{ 
(1-qY^{2})\,\eta_{R}(U_{n})-(u_{n}^{1/2}-u_{n}^{-1/2})-(u_{0}^{1/2}-u_{0}^{-1/2})q^{1/2}\,Y
\right\}, \\ 
\mathcal{R}_{i}(Y)&:=\frac{t^{1/2}}{1-tY}\left\{ 
(1-Y)\eta_{R}(T_{i})-(t^{1/2}-t^{-1/2}) \right\} \quad (1 \le i<n), \\ 
\mathcal{K}_{n}(Y)&:=\frac{t_{n}^{1/2}}{(1-t_{0}^{1/2}t_{n}^{1/2}Y^{-1})(1+t_{0}^{-1/2}t_{n}^{1/2}Y^{-1})} \\ 
&\times\left\{ 
(1-Y^{-2})\,\eta_{R}(T_{n})-(t_{n}^{1/2}-t_{n}^{-1/2})-(t_{0}^{1/2}-t_{0}^{-1/2})Y^{-1}
\right\}. 
\end{align*}
Then we have 
\begin{align}
C(1, \epsilon_{m})
&=\mathcal{R}_{m}(Y_{m+1}/Y_{m}) \cdots \mathcal{R}_{n-1}(Y_{n}/Y_{m})\mathcal{K}_{n}(Y_{m}) 
\label{eq:C-matrix} \\ 
&\times\mathcal{R}_{n-1}(Y_{m}^{-1}Y_{n}^{-1}) \cdots \mathcal{R}_{m}(Y_{m}^{-1}Y_{m+1}^{-1}) 
\nonumber \\ 
&\times\mathcal{R}_{m-1}(Y_{m-1}^{-1}Y_{m}^{-1}) \cdots \mathcal{R}_{1}(Y_{1}^{-1}Y_{m}^{-1}) 
\nonumber \\ 
&\times\mathcal{K}_{0}(Y_{m}^{-1})\mathcal{R}_{1}(qY_{1}/Y_{m}) \cdots \mathcal{R}_{m-1}(qY_{m-1}/Y_{m}). 
\nonumber
\end{align}




\subsection{The degenerate double affine Hecke algebra}

In this subsection we consider the trigonometric degeneration of the DAHA of type $(C_{n}^{\vee}, C_{n})$. 
We refer to \cite{C2} for the general theory on degeneration of the DAHA. 

In the following we make use of $X_{i}, T_{i}, Y_{i} \, (1\le i \le n)$ as generators of $\DAHA$.  
Note that $T_{0}$ is recovered from them by 
$T_{0}=T_{1}^{-1} \cdots T_{n-1}^{-1}\cdot T_{n}^{-1} \cdots T_{1}^{-1}Y_{1}$. 
Let $\hbar$ be a small parameter. 
We set 
\begin{align}
& 
q^{1/2}=e^{\hbar c/2}, \,\, t^{1/2}=e^{\hbar k/2}, 
\label{eq:trig-deg-1} \\ 
& 
t_{0}^{1/2}=e^{\hbar k_{0}/2}, \,\, 
t_{n}^{1/2}=e^{\hbar k_{n}/2}, \,\, u_{0}^{1/2}=e^{\hbar k_{0}^{*}/2}, \,\, 
u_{n}^{1/2}=e^{\hbar k_{n}^{*}/2}, 
\nonumber \\
& 
X_{i}=x_{i}, \quad T_{i}=s_{i}+\hbar\widetilde{T}_{i}+o(\hbar), \quad 
Y_{i}=e^{\hbar y_{i}} \qquad (1 \le i \le n).
\label{eq:trig-deg-2} 
\end{align}
and take the limit $\hbar \to 0$. 
In \eqref{eq:trig-deg-2} we introduced accessory generators $\widetilde{T}_{i} \, (1 \le i \le n)$ to avoid 
rewriting formulas in the form of Lusztig's relations. 
For example, substitute \eqref{eq:trig-deg-1} and \eqref{eq:trig-deg-2} into the relations  
\begin{align*}
(T_{i}-t^{1/2})(T_{i}+t^{-1/2})=0, \quad T_{i}Y_{i+1}T_{i}=Y_{i} \qquad (1 \le i \le n)  
\end{align*}
and expand them into power series of $\hbar$. 
Taking the zeroth and the first order terms we obtain 
\begin{align*}
s_{i}^{2}=1,  \quad s_{i}\widetilde{T}_{i}+\widetilde{T}_{i}s_{i}=ks_{i}, \quad 
\widetilde{T}_{i}s_{i}+s_{i}y_{i+1}s_{i}+s_{i}\widetilde{T}_{i}=y_{i}. 
\end{align*}
Eliminating $\widetilde{T}_{i}$, we find $y_{i}s_{i}=s_{i}y_{i+1}+k$. 
Thus we get closed relations among $x_{i}, s_{i}$ and $y_{i} \, (1 \le i \le n)$.  
The parameters $k_{0}, k_{n}, k_{0}^{*}$ and $k_{n}^{*}$ appear only in the form of  
$k_{0}+k_{n}$ and $k_{0}^{*}+k_{n}^{*}$. 
Setting $\alpha:=(k_{0}+k_{n})/2$ and $\beta:=(k_{0}^{*}+k_{n}^{*})/2$, 
we obtain the trigonometric degeneration of $\DAHA$: 

\begin{defn}\label{def:DAHA}
The {\it degenerate double affine Hecke algebra} $\dDAHA$ of type $(C_{n}^{\vee}, C_{n})$ 
is the unital associative algebra generated by $x_{i}, \, s_{i}$ and $y_{i} \,\, (1 \le i \le n)$ 
satisfying the following relations: 
\begin{align*}
& 
s_{i}^{2}=1 \quad (1 \le i \le n), \qquad 
s_{i}s_{i+1}s_{i}=s_{i+1}s_{i}s_{i+1} \quad (1 \le i <n), \\ 
& 
s_{n-1}s_{n}s_{n-1}s_{n}=s_{n}s_{n-1}s_{n}s_{n-1}, \\ 
& 
s_{i}x_{i}s_{i}=x_{i+1} \quad (1 \le i<n), \qquad  s_{n}x_{n}s_{n}=x_{n}^{-1}, \\ 
& 
y_{i}s_{i}=s_{i}y_{i+1}+k \quad (1 \le i<n), \qquad
y_{n}s_{n}=-s_{n}y_{n}+2\alpha, \\ 
& 
[s_{i}, x_{j}]=0=[s_{i}, y_{j}] \quad (|i-j|>1 \,\, \hbox{or} \,\, (i, j)=(n, n-1)),\\ 
& 
[y_{i}, x_{j}]=\left\{ 
\begin{array}{ll}
k(\tilde{s}_{ji}-s_{ji})x_{j} & (i>j) \\
cx_{i}+2(\alpha x_{i}^{-1}+\beta)r_{i} & \\ 
\quad\,\,{}+k(\sum_{1 \le l<i}s_{li}x_{l}+\sum_{i<l\le n}s_{il}x_{i}+\sum_{l (\not=i)}\tilde{s}_{il}x_{i}) & (i=j) \\ 
k(\tilde{s}_{ij}x_{j}-s_{ij}x_{i}) & (i<j), 
\end{array}
\right.
\end{align*}
where 
\begin{align*}
& 
s_{ij}=s_{ji}:=(s_{i} \cdots s_{j-1})(s_{j-2} \cdots s_{i}) \quad (i<j),  \\ 
& 
r_{i}:=(s_{i} \cdots s_{n})(s_{n-1} \cdots s_{i}) \quad (1 \le i \le n), \qquad 
\tilde{s}_{ij}:=r_{i}r_{j}s_{ij}. 
\end{align*}
\end{defn}

The subalgebra generated by $s_{i} \, (1 \le i <n)$ is isomorphic to 
the group algebra $\mathbb{C}[W_{0}]$. 
The subalgebra $\overline{H}$ generated by $s_{i} \, (1 \le i <n)$ and $y_{i} \, (1 \le i \le n)$ is called 
the {\it degenerate affine Hecke algebra} of type $C_{n}$. 
The subalgebra generated by $x_{i} \, (1 \le i \le n)$ and $s_{i} \, (1\le i<n)$ is 
isomorphic to the group algebra $\mathbb{C}[W]$ of the affine Weyl group
through the map $x_{1}^{-1}r_{1} \mapsto s_{0}$ and $s_{i} \mapsto s_{i} \, (1 \le i \le n)$. 
Hereafter we identify them. 
We have the PBW decomposition: 
\begin{align*}
\mathbb{C}[W] \simeq \mathbb{C}[x^{\pm 1}] \otimes \mathbb{C}[W_{0}], \quad
\overline{H} \simeq \mathbb{C}[W_{0}] \otimes \mathbb{C}[y], \quad 
\dDAHA \simeq 
\mathbb{C}[x^{\pm 1}] \otimes \mathbb{C}[W_{0}] \otimes 
\mathbb{C}[y],  
\end{align*}
where $\mathbb{C}[x^{\pm 1}]:=\mathbb{C}[x_{1}^{\pm 1}, \ldots , x_{n}^{\pm 1}]$ and 
$\mathbb{C}[y]:=\mathbb{C}[y_{1}, \ldots , y_{n}]$.

We regard $\overline{\mathbb{L}}:=\mathbb{C}[x^{\pm 1}] \otimes \mathbb{C}[y]$ 
as the ring of regular functions on 
$(\mathbb{C}^{\times})^{n} \times \mathbb{C}^{n}$.  
Denote by $\overline{\mathbb{K}}$ the field of meromorphic functions on $(\mathbb{C}^{\times})^{n} \times \mathbb{C}^{n}$.
We give $\overline{\mathbb{L}}$-module structure to $\dDAHA$ in the same way as \eqref{eq:L-mod-str}. 
Then $\mathbb{C}[W_{0}]^{\overline{\mathbb{K}}}:=\overline{\mathbb{K}}\otimes_{\overline{\mathbb{L}}}\dDAHA$ 
is the vector space of $\mathbb{C}[W_{0}]$-valued meromorphic functions. 
Any element $G \in \mathbb{C}[W_{0}]^{\overline{\mathbb{K}}}$ is uniquely represented 
in the form $G=\sum_{w \in W_{0}}g_{w}.w$ where $g_{w} \in \overline{\mathbb{K}}$. 

We define two maps 
\begin{align}
\overline{\eta}_{L}: \overline{H} \to 
{\rm End}_{\overline{\mathbb{K}}}(\mathbb{C}[W_{0}]^{\overline{\mathbb{K}}}), \qquad 
\overline{\eta}_{R}: \mathbb{C}[W] \to 
{\rm End}_{\overline{\mathbb{K}}}(\mathbb{C}[W_{0}]^{\overline{\mathbb{K}}})
\label{eq:bimod-action-eta}
\end{align}
by 
\begin{align*}
\overline{\eta}_{L}(h)\left(\sum_{w \in W_{0}}g_{w}.w \right):=
\sum_{w \in W_{0}}g_{w}.(h w), \quad
\overline{\eta}_{R}(\xi)\left(\sum_{w \in W_{0}}g_{w}.w \right):=
\sum_{w \in W_{0}}g_{w}.(w \xi).  
\end{align*}
Then the map $\overline{\eta}_{L}$ (resp.\,\,$\overline{\eta}_{R}$) 
is an algebra (resp.\,\,anti-algebra) homomorphism. 
Hence they determine
$(\overline{H}, \mathbb{C}[W])$-bimodule structure on $\mathbb{C}[W_{0}]^{\overline{\mathbb{K}}}$. 

\subsection{Degeneration of the bispectral qKZ equation ($Y$-side)}\label{subsec:y-side}

In the following two subsections, we consider the trigonometric degeneration of 
the bispectral qKZ equation \eqref{eq:QAKZ-X} and \eqref{eq:QAKZ-Y}. 

First we consider \eqref{eq:QAKZ-Y}.
Recall that the operator $C(1, \epsilon_{m})$ is explicitly given by \eqref{eq:C-matrix}. 
Denote by 
$\overline{\mathcal{K}}_{0}(y), \overline{\mathcal{R}}_{i}(y), \overline{\mathcal{K}}_{n}(y)$
the zeroth order term of the power series expansion of 
$\mathcal{K}_{0}(e^{\hbar y}), \mathcal{R}_{i}(e^{\hbar y}), \mathcal{K}_{n}(e^{\hbar y})$
as $\hbar \to 0$. 
We have 
\begin{align*}
\overline{\mathcal{K}}_{0}(y)&=\frac{1}{y+\beta+c/2}\left( 
\left(y+c/2\right)\overline{\eta}_{R}(x_{1}^{-1}r_{1})+\beta\right), \\ 
\overline{\mathcal{R}}_{i}(y)&=\frac{1}{y+k}\left(
y\,\overline{\eta}_{R}(s_{i})+k\right) \quad (1\le i<n), \\ 
\overline{\mathcal{K}}_{n}(y)&=\frac{1}{y-\alpha}\left( 
y\,\overline{\eta}_{R}(s_{n})-\alpha\right). 
\end{align*}

%

Suppose that $F=F(X\,|\,Y) \in H_{0}^{\mathbb{K}}$ is expanded as 
\begin{align}
F=\hbar^{M}\left( \overline{G}+o(1)\right) \qquad (\hbar \to 0)
\label{eq:solution-limit}
\end{align}
for some $M \in \mathbb{Z}$ and 
$\overline{G} \in \mathbb{C}[W_{0}]^{\overline{\mathbb{K}}}\setminus\{0\}$. 
Taking the lowest degree term with respect to $\hbar$, 
we obtain the trigonometric degeneration of \eqref{eq:QAKZ-Y}: 
\begin{align}
\overline{C}(1, \epsilon_{m})\Delta_{m}\overline{G}=\overline{G} \quad (1 \le m \le n),    
\label{eq:rational-qKZ-Hecke}
\end{align}
where 
\begin{align}
\overline{C}(1, \epsilon_{m})
&:=\overline{\mathcal{R}}_{m}(y_{m+1}-y_{m}) \cdots \overline{\mathcal{R}}_{n-1}(y_{n}-y_{m})
\overline{\mathcal{K}}_{n}(y_{m}) 
\label{eq:Cbar-matrix} \\ 
&\times\overline{\mathcal{R}}_{n-1}(-y_{m}-y_{n}) \cdots 
\overline{\mathcal{R}}_{m}(-y_{m}-y_{m+1}) 
\nonumber \\ 
&\times\overline{\mathcal{R}}_{m-1}(-y_{m-1}-y_{m}) \cdots 
\overline{\mathcal{R}}_{1}(-y_{1}-y_{m}) 
\nonumber \\ 
&\times\overline{\mathcal{K}}_{0}(-y_{m})
\overline{\mathcal{R}}_{1}(y_{1}-y_{m}+c) \cdots \overline{\mathcal{R}}_{m-1}(y_{m-1}-y_{m}+c). 
\nonumber
\end{align}

\subsection{Degeneration of the bispectral qKZ equation ($X$-side)}\label{subsec:x-side}

Next we consider the trigonometric degeneration of \eqref{eq:QAKZ-X}. 
One can calculate it by rewriting the operator $C(\epsilon_{a}, 1)$ 
into the form of \eqref{eq:C-matrix} and taking the limit as $\hbar \to 0$. 
Here we calculate the limit more directly. 

 
The QAKZ equation \eqref{eq:QAKZ-X} is equivalent to $\tau(\epsilon_{a}, 1)(F)=F \,\, (1 \le a \le n)$, 
which is explicitly given by 
\begin{align}
\sum_{w \in W_{0}}((\epsilon_{a}, 1)(f_{w}) \cdot d_{\epsilon_{a}}(X)^{-1}).\widetilde{S}_{\epsilon_{a}}T_{w}=
\sum_{w \in W_{0}}f_{w}.T_{w},
\label{eq:bispectral-qKZ-xside}
\end{align}
where $F=\sum_{w \in W_{0}}f_{w}.T_{w} \, (f_{w} \in \mathbb{K})$.
By direct calculation we have 
\begin{align*}
d_{\epsilon_{a}}(X)=d_{\epsilon_{a}}^{(0)}(x)
\left(1-\hbar\,
d_{\epsilon_{a}}^{\dagger}(x)+o(\hbar)\right), 
\end{align*}
where 
\begin{align*}
d_{\epsilon_{a}}^{(0)}(x):=(1-x_{a}^{2})^{2}
\prod_{p(\not=a)}(1-x_{a}/x_{p})(1-x_{a}x_{p})
\end{align*}
and
\begin{align*}
d_{\epsilon_{a}}^{\dagger}(x)&:=
\alpha\,\frac{1+x_{a}^{2}}{1-x_{a}^{2}}+2\beta\frac{x_{a}}{1-x_{a}^{2}} \\ 
&\,{}+c\left( 
\frac{x_{a}^{2}}{1-x_{a}^{2}}-\sum_{p=1}^{a-1}\frac{x_{a}}{x_{a}-x_{p}}
\right)+k\sum_{p(\not=a)}\left(
\frac{1}{1-x_{a}x_{p}}-\frac{x_{a}}{x_{a}-x_{p}}
\right).
\end{align*}
To calculate the limit of $\widetilde{S}_{\epsilon_{a}}$ as $\hbar \to 0$,  
we use the accessory generators $\widetilde{T}_{i} \,\, (1 \le i \le n)$ and the relations 
\begin{align*}
& 
s_{i}\widetilde{T}_{i}+\widetilde{T}_{i}s_{i}=ks_{i} \quad (1 \le i<n), \qquad 
s_{n}\widetilde{T}_{n}+\widetilde{T}_{n}s_{n}=k_{n}s_{n}, \\ 
& 
\widetilde{T}_{i}x_{i}=x_{i+1}(\widetilde{T}_{i}-k), \,\,
x_{i}\widetilde{T}_{i}=(\widetilde{T}_{i}-k)x_{i+1}  \quad (1 \le i<n), \\ 
& 
\widetilde{T}_{n}x_{n}=x_{n}^{-1}(\widetilde{T}_{i}-k_{n})-k_{n}^{*}  
\end{align*}
derived from 
(i) and (iii) in Definition \ref{def:DAHA}. 
Then we obtain 
\begin{align*}
\widetilde{S}_{\epsilon_{a}}=d_{\epsilon_{a}}^{(0)}(x)\left( 
1+\hbar\,\widetilde{S}_{\epsilon_{a}}^{\dagger}+o(\hbar)
\right), 
\end{align*}
where 
\begin{align*}
\widetilde{S}_{\epsilon_{a}}^{\dagger}&:=
\Phi_{a}-c\left(\frac{x_{a}^{2}}{1-x_{a}^{2}}-\sum_{p=1}^{a-1}\frac{x_{a}}{x_{a}-x_{p}}\right), \\
\Phi_{a}&:=y_{a}+2\,\frac{\alpha+\beta x_{a}}{x_{a}^{2}-1}r_{a} \\ 
&+k\left( 
\sum_{p=1}^{a-1}\frac{x_{a}}{x_{a}-x_{p}}s_{pa}+
\sum_{p=a+1}^{n}\frac{x_{p}}{x_{a}-x_{p}}s_{ap}+\sum_{p(\not=a)}
\frac{1}{x_{a}x_{p}-1}\tilde{s}_{ap}
\right). 
\end{align*}
Note that 
\begin{align*}
d_{\epsilon_{a}}^{\dagger}(x)+\widetilde{S}_{\epsilon_{a}}^{\dagger}=
\Phi_{a}+cx_{a}\frac{\partial \log{U(x)}}{\partial x_{a}}, 
\end{align*}
where 
\begin{align}
U(x)&:=\prod_{p=1}^{n}\left(x_{p}^{\frac{k(n-1)+\alpha}{c}}
(1+x_{p})^{\frac{\beta-\alpha}{c}}(1-x_{p})^{-\frac{\beta+\alpha}{c}}\right) 
\label{eq:gauge-U} \\ 
&\times 
\prod_{1\le p<l \le n}\bigl( 
(1-x_{p}x_{l})(x_{p}-x_{l})
\bigr)^{-\frac{k}{c}}.
\nonumber
\end{align}

Suppose that $F=F(X\,|\,Y) \in H_{0}^{\mathbb{K}}$ is expanded as in \eqref{eq:solution-limit}. 
Then we have 
\begin{align*}
(\epsilon_{a}, 1)F-F=\hbar^{M+1}(cx_{a}\frac{\partial \overline{G}}{\partial x_{a}}+o(1)). 
\end{align*}
Hence the trigonometric degeneration of \eqref{eq:bispectral-qKZ-xside} is given by 
\begin{align}
\left( cx_{a}\frac{\partial}{\partial x_{a}}+L_{a}+cx_{a}\frac{\partial \log{U(x)}}{\partial x_{a}}\right)
\overline{G}=0, 
\label{eq:differential-Hecke}
\end{align}
where 
\begin{align*}
L_{a}&:=\eta_{L}(y_{a})+2\,\frac{\alpha+\beta x_{a}}{x_{a}^{2}-1}\eta_{L}(r_{a}) \\ 
&+k\left( 
\sum_{p=1}^{a-1}\frac{x_{a}}{x_{a}-x_{p}}\eta_{L}(s_{pa})+
\sum_{p=a+1}^{n}\frac{x_{p}}{x_{a}-x_{p}}\eta_{L}(s_{ap})+
\sum_{p(\not=a)}\frac{1}{x_{a}x_{p}-1}\eta_{L}(\tilde{s}_{ap})
\right). 
\end{align*}

\begin{rem}
The equation \eqref{eq:differential-Hecke} is the semi-classical limit 
of the QAKZ equation \eqref{eq:QAKZ-X},  
and hence it could be regarded as 
the affine KZ equation of type $(C_{n}^{\vee}, C_{n})$. 
See Section 1.3.2 of \cite{C2} for details of 
such correspondence in the $GL_{n}$ case. 
\end{rem}

\subsection{Embedding into the compatible system}\label{subsec:embedding}

As was seen in subsections \ref{subsec:y-side} and \ref{subsec:x-side}, 
the trigonometric degeneration of the bispectral qKZ equation is the system of equations 
\eqref{eq:rational-qKZ-Hecke} and \eqref{eq:differential-Hecke} for 
$\overline{G} \in \mathbb{C}[W_{0}]^{\overline{\mathbb{K}}}$. 
Compatibility of the system formally follows from the relation 
$\tau(1, \epsilon_{m})(\tau(\epsilon_{a}, 1)-1)=(\tau(\epsilon_{a}, 1)-1)\tau(1, \epsilon_{m})$.   
In this subsection we prove that the system is contained in 
our compatible one \eqref{eq:main}. 
 
Define $G \in \mathbb{C}[W_{0}]^{\overline{\mathbb{K}}}$ by 
\begin{align}
G(x\,|\,y)=U(x)\overline{G}(x\,|\,y),  
\label{eq:gauge-transform}
\end{align}
where $U(x)$ is given by \eqref{eq:gauge-U}. 
Note that the operator $\overline{C}(1, \epsilon_{m})\Delta_{m}$ commutes with multiplication by
any function in $x$. 
Hence the modified function $G$ also satisfies \eqref{eq:rational-qKZ-Hecke}. 
Thus the system of equations \eqref{eq:rational-qKZ-Hecke} and \eqref{eq:differential-Hecke} 
is equivalent to 
\begin{align}
\overline{C}(1, \epsilon_{m})\Delta_{m}G
&=G \qquad (1 \le m \le n),  
\label{eq:difference-degeneration} \\
\left( cx_{a}\frac{\partial}{\partial x_{a}}+L_{a}\right)G
&=0 \qquad (1 \le a \le n)
\label{eq:differential-degeneration}
\end{align}
for $G \in \mathbb{C}[W_{0}]^{\overline{\mathbb{K}}}$. 


We realize the $(\overline{H}, \mathbb{C}[W])$-bimodule $\mathbb{C}[W_{0}]^{\overline{\mathbb{K}}}$ in 
the scalar extension 
$(V^{\otimes n})^{\overline{\mathbb{K}}}:=\overline{\mathbb{K}}\otimes_{\mathbb{C}}V^{\otimes n}$ 
as follows. 
First we define the left action of $W_{0}$ on $V$ by 
\begin{align*}
s_{i}(v_{a})=v_{\sigma_{i}(a)}, \quad 
s_{i}(v_{\overline{a}})=v_{\overline{\sigma_{i}(a)}} \qquad (1\le i<n), 
\end{align*}
where $\sigma_{i}:=(i, i+1)$ is the transposition, and 
\begin{align*}
s_{n}(v_{a})=\left\{
\begin{array}{ll}
v_{\overline{n}} & (a=n) \\
v_{a} & (a\not=n),
\end{array}
\right.
\qquad 
s_{n}(v_{\overline{a}})=\left\{
\begin{array}{ll}
v_{n} & (a=n) \\
v_{\overline{a}} & (a\not=n).
\end{array}
\right.
\end{align*}
View $V^{\otimes n}$ as a tensor representation of $W_{0}$. 
We extend it to the $\overline{\mathbb{K}}$-linear action on 
$(V^{\otimes n})^{\overline{\mathbb{K}}}$.  
Let $(V^{\otimes n})_{0}:=\mathbb{C}[W_{0}](v_{1}\otimes \cdots \otimes v_{n})$ 
be the cyclic submodule and 
$(V^{\otimes n})_{0}^{\overline{\mathbb{K}}}:=\overline{\mathbb{K}}\otimes_{\mathbb{C}}(V^{\otimes n})_{0}$.
Denote the left action on $(V^{\otimes n})_{0}^{\overline{\mathbb{K}}}$ by 
$\rho_{L}: \mathbb{C}[W_{0}] \to {\rm End}_{\overline{\mathbb{K}}}(V^{\otimes n})_{0}^{\overline{\mathbb{K}}}$.

Note that the operators $A_{a}(y), B_{a}(x)$ and $Q_{m}(x\,|\,y)$ are $\overline{\mathbb{K}}$-linear, 
and hence belong to ${\rm End}_{\overline{\mathbb{K}}}(V^{\otimes n})_{0}^{\overline{\mathbb{K}}}$. 
\begin{lem}\label{lem:AHA-action}
The following relations hold on $(V^{\otimes n})_{0}^{\overline{\mathbb{K}}}$: 
\begin{align*}
& 
A_{i}(y)\rho_{L}(s_{i})=\rho_{L}(s_{i})A_{i+1}(y)+k \quad (1\le i<n), \\
&
A_{n}(y)\rho_{L}(s_{n})=-\rho_{L}(s_{n})A_{n}(y)+2\alpha, \\ 
& 
[A_{i}(y), \rho_{L}(s_{j})]=0 \quad (|i-j|>1 \, \hbox{or} \, (i, j)=(n-1, n)).
\end{align*}
\end{lem}

{}From Lemma \ref{lem:AHA-action} the action $\rho_{L}$ is extended to 
that of $\overline{H}$ which maps $y_{a} \mapsto A_{a}(y)$. 
We also denote it by $\rho_{L}$. 
Now consider the $\overline{\mathbb{K}}$-linear map
\begin{align*}
\phi : \mathbb{C}[W_{0}]^{\overline{\mathbb{K}}} \to (V^{\otimes n})_{0}^{\overline{\mathbb{K}}}, \quad 
w \mapsto \rho_{L}(w)(v_{1}\otimes \cdots \otimes v_{n}) \quad (w \in W_{0}).
\end{align*}
 
\begin{prop}
The map $\phi$ gives an isomorphism between the left $\overline{H}$-modules 
$(\overline{\eta}_{L}, \mathbb{C}[W_{0}]^{\overline{\mathbb{K}}})$ and 
$(\rho_{L}, (V^{\otimes n})_{0}^{\overline{\mathbb{K}}})$.  
In particular we have $\phi \,\overline{\eta}_{L}(y_{a}) \phi^{-1}=A_{a}(y)$. 
\end{prop}

\begin{proof}
It is clear that $\phi$ commutes with the left action of $\mathbb{C}[W_{0}]$. 
{}From $A_{a}(y)(v_{1}\otimes \cdots \otimes v_{n})=y_{a} v_{1}\otimes \cdots \otimes v_{n}$ 
and Lemma \ref{lem:AHA-action}, $\phi$ commutes also with the action of 
$y_{a} \in \overline{H} \, (1 \le a \le n)$. 
\end{proof}

Recall that $\mathbb{C}[W_{0}]^{\overline{\mathbb{K}}}$ is a right $\mathbb{C}[W]$-module 
with the action $\overline{\eta}_{R}$ (see \eqref{eq:bimod-action-eta}). 
Now we define a right action of $\mathbb{C}[W]$ on $(V^{\otimes n})_{0}^{\overline{\mathbb{K}}}$: 

\begin{lem}\label{lem:right-action}
There exists an anti-algebra homomorphism 
$\rho_{R}:\mathbb{C}[W] \to {\rm End}_{\mathbb{\overline{K}}}(V^{\otimes n})_{0}^{\overline{\mathbb{K}}}$ 
such that 
\begin{align*}
\rho_{R}(s_{0})=T_{1}(x), \quad 
\rho_{R}(s_{i})=P_{i,i+1} \,\, (1\le i<n), \quad 
\rho_{R}(s_{n})=T_{n}(\underbar{$1$}).
\end{align*}  
\end{lem}

Lemma \ref{lem:right-action} follows from $T(x)^{2}=1$ and 
the braid relation $P_{i,i+1}P_{i+1,i+2}P_{i,i+1}=P_{i+1,i+2}P_{i,i+1}P_{i+1,i+2} \, (1 \le i \le n-2)$. 
We denote the right action of $\mathbb{C}[W]$ on $(V^{\otimes n})_{0}^{\overline{\mathbb{K}}}$ 
by $\rho_{R}$. 

\begin{prop}\label{prop:right-action}
The map $\phi$ commutes with the right action of $\mathbb{C}[W]$ 
on $\mathbb{C}[W_{0}]^{\overline{\mathbb{K}}}$ and 
$(V^{\otimes n})_{0}^{\overline{\mathbb{K}}}$. 
Therefore $\phi$ is an isomorphism between $(\overline{H}, \mathbb{C}[W])$-bimodules. 
\end{prop}

\begin{proof}
Set $v^{\dagger}:=v_{1} \otimes \cdots \otimes v_{n}$. 
For $w \in W_{0}$ we have 
\begin{align*}
\phi(w s_{i})=\rho_{L}(w)\rho_{L}(s_{i})v^{\dagger}=\rho_{L}(w)P_{i,i+1}v^{\dagger}=
P_{i,i+1}\rho_{L}(w)v^{\dagger}=P_{i,i+1}\phi(w)
\end{align*}
for $1\le i<n$, and 
\begin{align*} 
\phi(w s_{n})=\rho_{L}(w)\rho_{L}(s_{n})v^{\dagger}=\rho_{L}(w)T_{n}(\underbar{$1$})v^{\dagger}=
T_{n}(\underbar{$1$})\rho_{L}(w)v^{\dagger}=T_{n}(\underbar{$1$})\phi(w).     
\end{align*}
If $wx_{1}^{-1}$ is equal to $x_{j}^{-1}w$ (resp. $x_{j}w$) in $\mathbb{C}[W]$,  
the first component of $\rho_{L}(wr_{1})v^{\dagger} \in V^{\otimes n}$ is $v_{\bar{j}}$ (resp. $v_{j}$). 
Hence we get 
\begin{align*}
\phi (ws_{0})=\phi(w x_{1}^{-1}r_{1})=\phi(x_{j}^{-1}wr_{1})=x_{j}^{-1}\phi(wr_{1})=T_{1}(x)\rho_{L}(w)
v^{\dagger}=T_{1}(x)\phi(w).  
\end{align*}
\end{proof}

Now we send the equations \eqref{eq:difference-degeneration} and 
\eqref{eq:differential-degeneration} on $\mathbb{C}[W_{0}]$ by $\phi$. 
First consider \eqref{eq:difference-degeneration}. 
{}From \eqref{eq:Cbar-matrix} and Proposition \ref{prop:right-action} we find
\begin{align*}
& 
\phi \, \overline{C}(1, \epsilon_{m}) \phi^{-1} \\ 
&=p_{m}(y)^{-1}
\bigl((y_{m}-y_{m+1})P_{m,m+1}-k\bigr) \cdots 
\bigl((y_{m}-y_{n})P_{n-1,n}-k\bigr) \\ 
&\times 
\bigl(y_{m}T_{n}(\underbar{$1$})-\alpha\bigr)\bigl((y_{m}+y_{n})P_{n-1,n}-k\bigr) \cdots 
\bigl((y_{m}+y_{m+1})P_{m,m+1}-k\bigr) \\ 
&\times 
\bigl((y_{m}+y_{m-1})P_{m-1,m}-k\bigr) \cdots \bigl((y_{m}+y_{1})P_{1,2}-k\bigr)
\bigl((y_{m}-c/2)T_{1}(x)-\beta\bigr) \\ 
&\times 
\bigl((y_{m}-y_{1}-c)P_{1,2}-k\bigr) \cdots 
\bigl((y_{m}-y_{m-1}-c)P_{m-1,m}-k\bigr), 
\end{align*} 
where 
\begin{align*}
p_{m}(y)&:=(y_{m}-\alpha)(y_{m}-\beta-c/2) \\ 
&\times 
\prod_{j=1}^{m-1}(y_{m}-y_{j}-c-k)\prod_{j=m+1}^{n}(y_{m}-y_{j}-k)
\prod_{j(\not=m)}(y_{m}+y_{j}-k).  
\end{align*} 
Using 
\begin{align*}
\lambda P-k=(\lambda-k)PR(\lambda)^{-1}, \quad 
\lambda T(x)-\beta=(\lambda-\beta)K(\lambda \, | \, x, \beta)^{-1},  
\end{align*}
we obtain 
\begin{align*}
\phi \, \overline{C}(1, \epsilon_{m}) \phi^{-1}=Q_{m}(x\,|\,y)^{-1}. 
\end{align*}
Hence, sending \eqref{eq:difference-degeneration} by $\phi$, 
we get the boundary rational qKZ equation 
\begin{align*}
\Delta_{m}\phi(G)=Q_{m}(x\,|\,y)\phi(G) \qquad (1 \le m \le n)  
\end{align*}
on $(V^{\otimes n})_{0}$. 

Next let us consider \eqref{eq:differential-degeneration}. 
Note that the following equalities hold on $(V^{\otimes n})_{0}$:
\begin{align*}
& 
\sum_{j=1}^{n}\overline{E}_{aa}^{(j)}=\rho_{L}(r_{a}), \quad 
{\bf Y}_{aa}+{\bf Z}_{aa}=0, \\ 
& 
{\bf X}_{ab}+{\bf X}_{ba}=\rho_{L}(s_{ab}), \quad 
{\bf Y}_{ab}+{\bf Z}_{ab}=\rho_{L}(\tilde{s}_{ab}) \qquad (a\not=b). 
\end{align*}
Therefore we have 
\begin{align*}
\phi L_{a} \phi^{-1}=A_{a}(y)+B_{a}(x)=L_{a}(x\,|\,y).    
\end{align*}
The Euler operator $x_{a}\frac{\partial}{\partial x_{a}}$ commutes with $\phi$.  
Consequently, sending \eqref{eq:differential-degeneration} by $\phi$, 
we obtain the differential equation 
\begin{align*}
D_{a}(x \, | \, y)\phi(G)=0 \qquad (1 \le a \le n). 
\end{align*} 
As a result we find 

\begin{prop}
The trigonometric degeneration \eqref{eq:rational-qKZ-Hecke} and 
\eqref{eq:differential-Hecke} of the bispectral qKZ equation of type $(C_{n}^{\vee}, C_{n})$ 
is equivalent to the compatible system \eqref{eq:main} with $N=n$ restricted to $(V^{\otimes n})_{0}$ 
through the gauge transform \eqref{eq:gauge-transform}.
\end{prop}


\section{Integral formula for solutions in a special case}\label{sec:integral-formula}

We construct an integral formula for solutions to the compatible system \eqref{eq:main} 
with 
\begin{align*}
\alpha=\beta=k/2 
\end{align*}
and the variable $x$ restricted to the hyperplane 
\begin{align*}
x=(x_{1}, 1, \ldots , 1).  
\end{align*}
Hence only one differential operator $D_{1}(x \, | \, y)$ enters in the system. 
In the following we set 
\begin{align*}
x_{1}=e^{2\pi i \lambda} 
\end{align*}
and assume that 
\begin{align*}
\mathop{\rm Im}{c}>0, \qquad  
\mathop{\rm Im}{k}>0
\end{align*}
for simplicity's sake. 

Our solutions take values in the $2n$-dimensional subspace of $V^{\otimes n}$ 
determined as follows. 
Denote by $\tilde{V}$ the subspace of $V$ spanned by $v_{a}$ and $v_{\bar{a}} \, (2 \le a \le N)$. 
We fix a non-zero vector $\tilde{v} \in \tilde{V}^{\otimes (n-1)}$ satisfying 
\begin{align*}
P_{i,i+1}\tilde{v}=\tilde{v} \quad (1\le \forall{i} \le n-2), \qquad 
T_{n-1}(\underbar{$1$})\,\tilde{v}=\tilde{v}.  
\end{align*}
Now define the vectors $u_{j} \in V^{\otimes n} \, (1 \le j \le 2n)$ by 
\begin{align*}
u_{j}:=P_{1,j}(v_{1} \otimes \tilde{v}), \quad 
u_{2n+1-j}:=P_{j,n}(\tilde{v} \otimes v_{\bar{1}}) \qquad (1 \le j \le n),  
\end{align*}
where $P_{1,1}=P_{n,n}=id$. 
Our solutions take values in the subspace $\mathcal{V}:=\oplus_{j=1}^{2n}\mathbb{C}u_{j}$. 

Define the rational functions $g_{j}(t)=g_{j}(t\,|y) \,\, (1 \le j \le 2n)$ by 
\begin{align*}
g_{j}(t)&:=\frac{1}{t-y_{j}}\prod_{p=1}^{j-1}\frac{t-y_{p}-k}{t-y_{p}}, \\ 
g_{2n+1-j}(t)&:=\frac{1}{t+y_{j}}\prod_{p=j+1}^{n}\frac{t+y_{p}-k}{t+y_{p}}\,
\prod_{p=1}^{n}\frac{t-y_{p}-k}{t-y_{p}},
\end{align*}
for $1 \le j \le n$. 
Set $\mathcal{G}:=\sum_{j=1}^{2n}{\mathbb{C}}g_{j}(t)$.  

Denote by $\mathcal{W}$ the $\mathbb{C}$-vector space spanned by the functions in the form 
\begin{align}
\frac{P(e^{2\pi it/c})}{\prod_{p=1}^{n}(1-e^{2\pi i(t-y_{p})/c})(1-e^{2\pi i(t+y_{p})/c})}, 
\label{eq:deformed-cycle}
\end{align}
where $P$ is a polynomial whose coefficients are entire and periodic functions 
in $y_{1}, \ldots , y_{n}$ with period $c$. 

Now we define a pairing $I$ between $\mathcal{G}$ and $\mathcal{W}$ by 
\begin{align}
I(g, W):=\int_{C(y)}\varphi(t\,|\,y)g(t)W(e^{2\pi i t/c})\,dt \qquad 
(g \in \mathcal{G}, \, W \in \mathcal{W}), 
\label{eq:pairing}
\end{align}
where the kernel function $\varphi$ is defined by 
\begin{align*}
\varphi(t\,|\,y):=e^{-\frac{2\pi i \lambda}{c}t}
\prod_{p=1}^{n}
\frac{\Gamma\left(\frac{t-y_{p}-k}{-c}\right)\Gamma\left(\frac{t+y_{p}-k}{-c}\right)}
     {\Gamma\left(\frac{t-y_{p}}{-c}\right)\Gamma\left(\frac{t+y_{p}}{-c}\right)}. 
\end{align*}
The contour $C(y)$ is a deformation of the real line $(-\infty, +\infty)$ such that 
the poles at $\pm y_{p}+k+c\mathbb{Z}_{\ge 0} \,\, (1 \le p \le n)$ 
are above $C(y)$ and the poles at $\pm y_{p}+c\mathbb{Z}_{\le 0} \,\, (1 \le p \le n)$ are below $C(y)$. 
Suppose that $W(e^{2\pi i t/c})$ is given by \eqref{eq:deformed-cycle}. 
{}From the Stirling formula we see that the integral \eqref{eq:pairing} converges if 
\begin{align}
\mathop{\rm Re}{\lambda}<\deg{P}<\mathop{\rm Re}{\lambda}+2n. 
\label{eq:degree}
\end{align}

For $W \in \mathcal{W}$ satisfying the degree condition \eqref{eq:degree} we set 
\begin{align}
f(\lambda \, | \, y):=\sum_{j=1}^{2n}I(g_{j}, W)u_{j}.  
\label{eq:integral-formula}
\end{align}

\begin{prop}\label{prop:integral-1}
The function $f$ satisfies the boundary rational qKZ equation \eqref{eq:boundary-rational-qKZ} 
with $\alpha=\beta=k/2$ and $x=(e^{2\pi i\lambda}, 1, \ldots , 1)$. 
\end{prop}

\begin{proof}
In the proof below we need to signify the dependence of $W \in \mathcal{W}$ on $y$. 
For that purpose we set 
\begin{align*}
& 
\tilde{g}(t\,|\,y')=\sum_{j=1}^{2n}g_{j}(t\,|\,y')u_{j}, \\ 
& 
\tilde{f}(\lambda \,|\,y'\,|\,y)=
\int_{C(y')}\varphi(t\,|\,y')\tilde{g}(t\,|\,y')W(e^{2\pi i t/c}\,|\,y)dt, 
\end{align*}
for $y'=(y_{1}', \ldots , y_{n}')$ such that each coordinate $y_{j}'$ belongs to 
the set $\{ \pm y_{a}+cl \, |\, 1 \le a \le n, \, l \in \mathbb{Z}\}$. 

By direct calculation we get 
\begin{align*}
P_{l,l+1}R_{l,l+1}(y_{l}-y_{l+1})\, \tilde{g}(t\,|\, \ldots , y_{l}, y_{l+1}, \ldots)=
\tilde{g}(t\,|\, \ldots , y_{l+1}, y_{l}, \ldots) 
\end{align*}
for $1 \le l \le n-1$, and 
\begin{align*} 
K_{n}(y_{n} \,|\, \underbar{$1$}, k/2)\,\tilde{g}(t\,|\, y_{1}, \ldots , y_{n-1}, y_{n})=
\tilde{g}(t\,|\, y_{1}, \ldots , y_{n-1}, -y_{n}).  
\end{align*}
Since $\varphi(t\,|\,y)$  and $C(y)$ are invariant under 
the transposition $y_{l} \leftrightarrow y_{l+1}$ and the reflection $y_{n} \to -y_{n}$, 
the above equalities where $\tilde{g}(t\,|\,\ldots)$ is replaced by 
$\tilde{f}(\lambda\,| \ldots |\, y)$ also hold. 
{}From this fact and the periodicity of $W$ in $y$, 
it is enough to prove that 
\begin{align}
K_{1}(y_{1}-c/2\,|\, x, k/2)\tilde{f}(\lambda\,|\,-y_{1}, y_{2}, \ldots , y_{n}\,|\, y)=
\tilde{f}(\lambda\,|\,y_{1}-c, y_{2}, \ldots , y_{n} \,|\, y) 
\label{eq:qkz-proof}
\end{align}
with $x=(e^{2\pi i \lambda}, 1, \ldots , 1)$. 
We abbreviate $y^{(1)}=(-y_{1}, y_{2}, \ldots , y_{n})$ and 
$y^{(2)}=(y_{1}-c, y_{2}, \ldots , y_{n})$. 
Taking the coefficients of $u_{j} \, (1 \le j \le 2n)$ in \eqref{eq:qkz-proof} we have the equalities to prove: 
\begin{align}
& \label{eq:qkz-proof-1}
\int_{C(y^{(1)})}\varphi(t\,|\,y^{(1)})g_{j}(t\,|\,y^{(1)})W(e^{2\pi i t/c})dt \\ 
& \nonumber {}=
\int_{C(y^{(2)})}\varphi(t\,|\,y^{(2)})g_{j}(t\,|\,y^{(2)})W(e^{2\pi i t/c})dt \qquad (j\not=1, 2n), \\ 
& \label{eq:qkz-proof-2}
\int_{C(y^{(1)})}\varphi(t\,|\,y^{(1)})\,
\frac{e^{2\pi i \lambda}(y_{1}-c/2)g_{2n}(t\,|\,y^{(1)})+(k/2)g_{1}(t\,|\,y^{(1)})}{y_{1}-c/2+k/2}\,
W(e^{2\pi i t/c})dt \\
& \nonumber {}=
\int_{C(y^{(2)})}\varphi(t\,|\,y^{(2)})g_{1}(t\,|\,y^{(2)})W(e^{2\pi i t/c})dt, \\ 
& \label{eq:qkz-proof-3}
\int_{C(y^{(1)})}\varphi(t\,|\,y^{(1)})\,
\frac{e^{-2\pi i \lambda}(y_{1}-c/2)g_{1}(t\,|\,y^{(1)})+(k/2)g_{2n}(t\,|\,y^{(1)})}{y_{1}-c/2+k/2}\,
W(e^{2\pi i t/c})dt \\
& \nonumber {}=
\int_{C(y^{(2)})}\varphi(t\,|\,y^{(2)})g_{2n}(t\,|\,y^{(2)})W(e^{2\pi i t/c})dt. 
\end{align}

First we prove \eqref{eq:qkz-proof-1}. 
If $j \not=1, 2n$, we have 
\begin{align}
\frac{g_{j}(t\,|\,y^{(1)})}{g_{j}(t\,|\,y^{(2)})}=\frac{t+y_{1}-k}{t+y_{1}}
\frac{t-y_{1}+c}{t-y_{1}-k+c}.
\label{eq:g-fraction}
\end{align}
Note that $\varphi(t\,|\,y^{(1)})=\varphi(t\,|\,y)$ and  
\begin{align}
\frac{\varphi(t\,|\,y^{(2)})}{\varphi(t\,|\,y)}=\frac{t+y_{1}-k}{t+y_{1}}\frac{t-y_{1}+c}{t-y_{1}-k+c}.  
\label{eq:qkz-proof-shift-y}
\end{align}
It is equal to the right hand side of \eqref{eq:g-fraction}, and hence
\begin{align*}
\varphi(t\,|\,y^{(1)})g_{j}(t\,|\,y^{(1)})=\varphi(t\,|\,y^{(2)})g_{j}(t\,|\,y^{(2)}) 
\qquad (j\not=1, 2n). 
\end{align*}
Thus the integrands in the both hand sides of \eqref{eq:qkz-proof-1} are the same. 
Since the integrand has no poles at $-y_{1}+k, -y_{1}+c, y_{1}$ and $y_{1}+k-c$, 
we can deform $C(y^{(1)})$ to $C(y^{(2)})$ without crossing poles. 
Thus we obtain \eqref{eq:qkz-proof-1}. 


In the rest we only prove \eqref{eq:qkz-proof-2}. 
The proof of \eqref{eq:qkz-proof-3} is similar. 
We separate the integrand in the left hand side and first consider 
\begin{align}
\frac{y_{1}-c/2}{y_{1}-c/2+k/2}
\int_{C(y^{(1)})}\varphi(t\,|\,y^{(1)})\,e^{2\pi i \lambda}g_{2n}(t\,|\,y^{(1)})W(e^{2\pi i t/c})dt. 
\label{eq:qkz-proof-tmp1}
\end{align}
We have 
\begin{align*}
\varphi(t\,|\,y^{(1)})\,e^{2\pi i \lambda}g_{2n}(t\,|\,y^{(1)})=\varphi(t-c\,|\,y)\,\frac{1}{t-y_{1}-k}.
\end{align*}
Changing $t \to t+c$ we see that the integral \eqref{eq:qkz-proof-tmp1} is equal to  
\begin{align*}
\frac{y_{1}-c/2}{y_{1}-c/2+k/2}\int_{C(y^{(1)})-c}
\varphi(t\,|\,y)\,\frac{1}{t+c-y_{1}-k}\,W(e^{2\pi i t/c})dt.  
\end{align*}
Since the integrand has no poles at $\pm y_{j} \, (1 \le j \le n)$, the contour $C(y^{(1)})-c$ 
can be deformed to the contour $C'$, 
which is a deformation of the real line such that 
the points $-y_{1}+k+c\mathbb{Z}_{\ge 0}, y_{1}+k+c\mathbb{Z}_{\ge -1}$ and 
$\pm y_{j}+k+c\mathbb{Z}_{\ge 0} \, (2 \le j \le n)$ are above $C'$, 
and the points $-y_{1}+c\mathbb{Z}_{\le 0}, y_{1}+c\mathbb{Z}_{\le -1}$ and 
$\pm y_{j}+c\mathbb{Z}_{\le 0} \, (2 \le j \le n)$ are below $C'$.  
The integrand of the rest part of the left hand side of \eqref{eq:qkz-proof-2} has no pole at 
$y_{1}$, hence the contour $C(y^{(1)})$ can be deformed to $C'$. 
As a result the left hand side becomes 
\begin{align*}
& 
\int_{C'}\varphi(t\,|\,y)\left(
\frac{y_{1}-c/2}{y_{1}-c/2+k/2}\frac{1}{t+c-y_{1}-k}+\frac{k/2}{y_{1}-c/2+k/2}\frac{1}{t+y_{1}}\right)
W(e^{2\pi i t/c})dt \\ 
& {}=
\int_{C'}\varphi(t\,|\,y)\frac{1}{t+c-y_{1}-k}\frac{t+y_{1}-k}{t+y_{1}}W(e^{2\pi i t/c})dt.
\end{align*}
On the other hand, using \eqref{eq:qkz-proof-shift-y} 
we see that the right hand side of \eqref{eq:qkz-proof-2} is equal to 
\begin{align*}
\int_{C(y^{(2)})}\varphi(t\,|\,y)\,\frac{1}{t+c-y_{1}-k}\frac{t+y_{1}-k}{t+y_{1}}\,W(e^{2\pi i t/c})dt. 
\end{align*} 
Since the integrand is regular at $t=-y_{1}+c$, 
the contour $C(y^{(2)})$ can be deformed to $C'$ without crossing poles. 
Thus we get the equality \eqref{eq:qkz-proof-2}. 
\end{proof}

\begin{prop}\label{prop:integral-2}
If $\alpha=\beta=k/2$ and $x=(e^{2\pi i \lambda}, 1, \ldots , 1)$, we have 
\begin{align*}
D_{1}(x\,|\,y)f(\lambda\,|\,y)=-\frac{ke^{2\pi i \lambda}}{e^{2\pi i \lambda}-1}f(\lambda\,|\,y).  
\end{align*}  
\end{prop}

\begin{proof}
Since ${\bf Y}_{11}$ and ${\bf Z}_{11}$ act as zero on $\mathcal{V}$, we find 
\begin{align*}
D_{1}(x\,|\,y)|_{\mathcal{V}}=\frac{c}{2\pi i}\frac{\partial}{\partial \lambda}+
\sum_{j=1}^{n}y_{j}(e_{11}^{(j)}-e_{\bar{1}\bar{1}}^{(j)}) 
&+\frac{k}{e^{2\pi i \lambda}-1}\left(
\sum_{j=1}^{n}e_{\bar{1}1}^{(j)}+\sum_{p=2}^{N}({\bf X}_{p1}+{\bf Z}_{1p})\right) \\ 
&+\frac{ke^{2\pi i \lambda}}{e^{2\pi i \lambda}-1}\left(\sum_{j=1}^{n}e_{1\bar{1}}^{(j)}+
\sum_{p=2}^{N}({\bf X}_{1p}+{\bf Y}_{1p})\right) 
\end{align*}
if $x=(e^{2\pi i \lambda}, 1, \ldots , 1)$ and $\alpha=\beta=k/2$. 
For $1 \le j \le n$ we have 
\begin{align*}
\sum_{p=2}^{N}({\bf X}_{p1}+{\bf Z}_{1p})v_{j}=\sum_{l=j+1}^{2n-j}v_{l}+\sum_{l=2n+2-j}^{2n}\!\!\!v_{l}, 
\quad 
\sum_{p=2}^{N}({\bf X}_{1p}+{\bf Y}_{1p})v_{j}=\sum_{l=1}^{j-1}v_{l}
\end{align*}
and 
\begin{align*}
\sum_{p=2}^{N}({\bf X}_{p1}+{\bf Z}_{1p})v_{2n+1-j}=\sum_{l=2n+2-j}^{2n}\!\!\!v_{l}, \quad 
\sum_{p=2}^{N}({\bf X}_{1p}+{\bf Y}_{1p})v_{2n+1-j}=\sum_{l=1}^{j-1}v_{l}+\sum_{l=j+1}^{2n-j}v_{l}. 
\end{align*}
{}From the calculation above we obtain 
\begin{align*}
D_{1}(x\,|\,y)f(\lambda\,|\,y)=\sum_{j=1}^{2n}I(h_{j}, W)u_{j}, 
\end{align*} 
where the rational functions $h_{j} \, (1\le j \le 2n)$ are given by 
\begin{align*}
& 
h_{j}(t\,|\,y)=-(t-y_{j})g_{j}(t\,|\,y)+
\frac{k}{e^{2\pi i \lambda}-1}\sum_{l=1}^{j-1}g_{l}(t\,|\,y)+
\frac{ke^{2\pi i \lambda}}{e^{2\pi i \lambda}-1}\sum_{l=j+1}^{2n}g_{l}(t\,|\,y), \\ 
& 
h_{2n+1-j}(t\,|\,y)=-(t+y_{j})g_{j}(t\,|\,y)+
\frac{k}{e^{2\pi i \lambda}-1}\sum_{l=1}^{2n-j}g_{l}(t\,|\,y)+
\frac{ke^{2\pi i \lambda}}{e^{2\pi i \lambda}-1}\sum_{l=2n+2-j}^{2n}\!\!\!g_{l}(t\,|\,y) 
\end{align*}
for $1 \le j \le n$. 

Note that 
\begin{align*}
g_{j}(t\,|\,y)&=\frac{1}{k}\left( \prod_{l=1}^{j-1}\frac{t-y_{l}-k}{t-y_{l}}-\prod_{l=1}^{j}\frac{t-y_{l}-k}{t-y_{l}}\right), \\ 
g_{2n+1-j}(t\,|\,y)&=\frac{1}{k}\left( 
\prod_{l=j+1}^{n}\frac{t+y_{l}-k}{t+y_{l}}-\prod_{l=j}^{n}\frac{t+y_{l}-k}{t+y_{l}} \right)
\prod_{l=1}^{n}\frac{t-y_{l}-k}{t-y_{l}}
\end{align*}
for $1 \le j \le n$. 
Using these formulas we get 
\begin{align*}
h_{j}(t\,|\,y)=-\frac{ke^{2\pi i \lambda}}{e^{2\pi i \lambda}-1}g_{j}(t\,|\,y)+\frac{1}{1-e^{2\pi i \lambda}}\left( 
1-e^{2\pi i \lambda}\prod_{l=1}^{n}\frac{t-y_{l}-k}{t-y_{l}}\frac{t+y_{l}-k}{t+y_{l}}
\right) 
\end{align*}
for $1 \le j \le 2n$. 
The last term in the right hand side is related to the kernel function $\varphi$ by 
\begin{align*}
e^{2\pi i \lambda}\prod_{l=1}^{n}\frac{t-y_{l}-k}{t-y_{l}}\frac{t+y_{l}-k}{t+y_{l}}=
\frac{\varphi(t-c\,|\,y)}{\varphi(t\,|\,y)}.  
\end{align*}
Therefore we have 
\begin{align}
& 
\int_{C(y)}\varphi(t\,|\,y)\left( 
1-e^{2\pi i \lambda}\prod_{l=1}^{n}\frac{t-y_{l}-k}{t-y_{l}}\frac{t+y_{l}-k}{t+y_{l}}
\right)W(e^{2\pi i t/c})dt  \label{eq:diff-proof-1} \\ 
&=
\int_{C(y)}\left(\varphi(t\,|\,y)-\varphi(t-c\,|\,y)\right)W(e^{2\pi it/c})dt 
\nonumber \\ 
&= 
\left(\int_{C(y)}-\int_{C(y)-c}\right)\varphi(t\,|\,y)W(e^{2\pi i t/c})dt. 
\nonumber
\end{align}
Since the integrand $\varphi(t\,|\,y)W(e^{2\pi i t/c})$ is regular at $t=\pm y_{j} \, (1 \le j \le n)$, 
the integral \eqref{eq:diff-proof-1} is equal to zero. 
Therefore we obtain 
\begin{align*}
D_{1}(x\,|\,y)f(e^{2\pi i \lambda}\,|\,y)&=\sum_{j=1}^{2n}I(h_{j}, W)u_{j} \\ 
&=-\frac{ke^{2\pi i \lambda}}{e^{2\pi i \lambda}-1}
\sum_{j=1}^{2n}I(g_{j}, W)u_{j}=-\frac{ke^{2\pi i \lambda}}{e^{2\pi i \lambda}-1}f(\lambda\,|\,y). 
\end{align*} 
\end{proof}

The linear operator $L_{1}(x\,|\,y)$ commutes with multiplication by any function in $x$. 
Hence, from Proposition \ref{prop:integral-1} and Proposition \ref{prop:integral-2}, 
we finally get 
\begin{thm}
Set 
\begin{align*}
\tilde{f}(\lambda\,|\,y):=(e^{2\pi i \lambda}-1)^{k/c}f(\lambda\,|\,y), 
\end{align*}
where $f$ is defined by \eqref{eq:integral-formula}. 
Then the function $\tilde{f}$ is a solution to the compatible system \eqref{eq:main} 
with $\alpha=\beta=k/2$ and $x=(e^{2\pi i \lambda}, 1, \ldots , 1)$ for any $W \in \mathcal{W}$ 
satisfying \eqref{eq:degree}. 
\end{thm}


\bigskip 
\noindent{\it Acknowledgements.} 
The research of the author is supported by Grant-in-Aid for 
Young Scientists (B) No.\,20740088. 
The author is deeply grateful to Saburo Kakei, Masahiro Kasatani, Michitomo Nishizawa and Yoshihisa Saito 
for valuable discussions. 
He also thanks the referee for critical comments that improved the paper considerably.

\end{document}